\documentclass[reqno,oneside,12pt]{amsart}

\hoffset - 1.5cm

\usepackage{amsmath}
\usepackage{amssymb}
\usepackage{amsthm}
\usepackage{amstext}
\usepackage{amsopn}
\usepackage{amsfonts}

\newtheorem{Fact}{Fact}[section]
\newtheorem{Lemma}{Lemma}[section]
\newtheorem{Theorem}{Theorem}[section]
\newtheorem{Proposition}{Proposition}[section]
\newtheorem{Hipothesis}{Hipothesis}[section]

\theoremstyle{definition}
\newtheorem{Corollary}{Corollary}[section]
\newtheorem{Definition}{Definition}[section]
\newtheorem{Example}{Example}[section]
\newtheorem{Remark}{Remark}[section]

\newcommand{\ba}{\begin{array}}
\newcommand{\bc}{\begin{center}}
\newcommand{\bd}{\begin{description}}
\newcommand{\bdm}{\begin{displaymath}}
\newcommand{\be}{\begin{enumerate}}
\newcommand{\beq}{\begin{equation}}
\newcommand{\bdf}{\begin{Definition}}
\newcommand{\bex}{\begin{Example}}
\newcommand{\bft}{\begin{Fact}}
\newcommand{\bl}{\begin{Lemma}}
\newcommand{\bp}{\begin{Proposition}}
\newcommand{\br}{\begin{Remark}}
\newcommand{\bt}{\begin{Theorem}}
\newcommand{\bco}{\begin{Corollary}}
\newcommand{\bh}{\begin{Hipothesis}}
\newcommand{\ea}{\end{array}}
\newcommand{\ec}{\end{center}}
\newcommand{\ed}{\end{description}}
\newcommand{\edm}{\end{displaymath}}
\newcommand{\ee}{\end{enumerate}}
\newcommand{\eeq}{\end{equation}}
\newcommand{\edf}{\end{Definition}}
\newcommand{\eex}{\end{Example}}
\newcommand{\eft}{\end{Fact}}
\newcommand{\el}{\end{Lemma}}
\newcommand{\ep}{\end{Proposition}}
\newcommand{\er}{\end{Remark}}
\newcommand{\et}{\end{Theorem}}
\newcommand{\eco}{\end{Corollary}}
\newcommand{\eh}{\end{Hipothesis}}

\newcommand{\bH}{\mathbb{H}}
\newcommand{\bI}{\mathbb{I}}

\newcommand{\bN}{\mathbb{N}}

\newcommand{\bR}{\mathbb{R}}

\newcommand{\bV}{\mathbb{V}}
\newcommand{\bW}{\mathbb{W}}

\newcommand{\bZ}{\mathbb{Z}}

\newcommand{\cC}{\mathcal{C}}

\newcommand{\cK}{\mathcal{K}}

\newcommand{\cU}{\mathcal{U}}

\newcommand{\numsec}{\setcounter{Theorem}{0}\setcounter{Definition}{0}
\setcounter{Remark}{0} \setcounter{Lemma}{0} \setcounter{Fact}{0}
\setcounter{Proposition}{0} \setcounter{Corollary}{0}
\setcounter{Example}{0} \setcounter{equation}{0}
\setcounter{Property}{0}\renewcommand\theequation{\arabic{section}.\arabic{equation}}
\renewcommand\theTheorem{\arabic{section}.\arabic{Theorem}}
\renewcommand\theDefinition{\arabic{section}.\arabic{Definition}}
\renewcommand\theRemark{\arabic{section}.\arabic{Remark}}
\renewcommand\theLemma{\arabic{section}.\arabic{Lemma}}
\renewcommand\theFact{\arabic{section}.\arabic{Fact}}
\renewcommand\theProposition{\arabic{section}.\arabic{Proposition}}
\renewcommand\theCorollary{\arabic{section}.\arabic{Corollary}}
\renewcommand\theExample{\arabic{section}.\arabic{Example}}
\renewcommand\theProperty{\arabic{section}.\arabic{Property}}}

\newcommand{\numsubsec}{\setcounter{Theorem}{0}\setcounter{Definition}{0}
\setcounter{Remark}{0} \setcounter{Lemma}{0} \setcounter{Fact}{0}
\setcounter{Proposition}{0} \setcounter{Corollary}{0}
\setcounter{Example}{0} \setcounter{equation}{0}
\setcounter{Property}{0}\renewcommand\theequation{\arabic{section}.\arabic{subsection}.\arabic{equation}}
\renewcommand\theTheorem{\arabic{section}.\arabic{subsection}.\arabic{Theorem}}
\renewcommand\theDefinition{\arabic{section}.\arabic{subsection}.\arabic{Definition}}
\renewcommand\theRemark{\arabic{section}.\arabic{subsection}.\arabic{Remark}}
\renewcommand\theLemma{\arabic{section}.\arabic{subsection}.\arabic{Lemma}}
\renewcommand\theFact{\arabic{section}.\arabic{subsection}.\arabic{Fact}}
\renewcommand\theProposition{\arabic{section}.\arabic{subsection}.\arabic{Proposition}}
\renewcommand\theCorollary{\arabic{section}.\arabic{subsection}.\arabic{Corollary}}
\renewcommand\theExample{\arabic{section}.\arabic{subsection}.\arabic{Example}}
\renewcommand\theProperty{\arabic{section}.\arabic{subsection}.\arabic{Property}}}

\numberwithin{equation}{section} \errorcontextlines=0
\newcommand{\im}{\mathrm{ im \;}}
\newcommand{\diag}{\mathrm{ diag \;}}

\newcommand{\sign}{\mathrm{ sign \;}}

\newcommand{\sone}{SO(2)}
\newcommand{\ds}{\displaystyle}
\newcommand{\nt}{\noindent}

\newcommand{\h}{\mathbb{H}}
\newcommand{\dg}{\nabla_{\sone}\mathrm{-deg}}
\newcommand{\dls}{\mathrm{deg_{LS}}}

\def\ra{\rightarrow}
\def\rg{\rangle}
\def\lg{\langle}
\def\ep{\epsilon}
\def\p{|}
\def\VS{\bV_{-\Delta}}
\def\nn{\nonumber}

\def\h1{{H^1(\o)}}
\def\hh1{{\mathbb{H}^1(\o)}}
\def\o{\Omega}

\def\si{{\sigma(-\Delta;\o)}}

\newcommand{\minev}{\lambda_0}
\def\bif{\textsc{Bif}(\infty,[\lambda_-,\lambda_+])}

\newenvironment{Proof}{\begin{proof}}{\end{proof}}

\begin{document}

\title[Neumann problem]{Existence and continuation of solutions \\ for a nonlinear Neumann problem}

\author{Krzysztof Muchewicz$^{\dag}$}
\address{
Faculty of Mathematics and Computer Science\\
Nicolaus Copernicus University \\
PL-87-100 Toru\'{n} \\ ul. Chopina 12/18 \\
Poland} \email{Krzysztof.Muchewicz@mat.uni.torun.pl}

\author{S{\l}awomir Rybicki$^{\ddag}$}
\address{Faculty of Mathematics and Computer Science\\
Nicolaus Copernicus University \\
PL-87-100 Toru\'{n} \\ ul. Chopina 12/18 \\
Poland} \email{Slawomir.Rybicki@mat.uni.torun.pl}

\date{\today}

\thanks{$^{\dag}$ Research sponsored by the Doctoral Program in Mathematics at the Nicolaus Copernicus University,
Toru\'n, Poland}

\thanks{$^{\ddag}$Partially supported by the Ministry Education and Science, Poland,  under grant 1 PO3A 009 27}


\begin{abstract}
In this article we study the existence, continuation  and bifurcation from infinity  of nonconstant solutions for a
nonlinear Neumann problem. We apply the Leray-Schauder degree and   the   degree for $\sone$-equivariant gradient
operators defined by the second author in \cite{[RYB1]}.
\end{abstract}
\maketitle


\numsec
\section{Introduction}

Consider the following nonlinear Neumann problem
\begin{equation}\label{intro:eq:1}
\left\{\begin{array}{rcll} -\Delta u  &=& f(u)& \mathrm{in} \;\o, \\
\ds \frac{\partial u}{\partial \nu} &=& 0 & \mathrm{on} \;
\partial \o,
\end{array}\right.
\end{equation}
where $\o \subset \bR^n$ is  an open bounded domain with $C^{1_-}$-boundary    and $f\in C^1(\bR,\bR)$.

The existence and multiplicity of weak solutions of problem \eqref{intro:eq:1} has been studied by many authors, see
for instance Hirano and Wan \cite{[HIRANOWAN]}, Ko \cite{[KO]}, Li \cite{[LI]}, Li and Li \cite{[LILI]}, Pomponio
\cite{[POMP]}, Tang \cite{[TANG]}, Tang and Wu  \cite{[TANGWU0]}, \cite{[TANGWU]}, Vanella\cite{[VANNELLA0]} and
references therein.

Usually weak solutions of system \eqref{intro:eq:1}  are considered as critical points of  a   functional $\Phi \in
C^2(\h1,\bR).$ The authors     apply tools of the critical point theory, like the Morse theory, the Conley index
technique and the mountain pass theorem, to obtain   results.

Solutions  of problem \eqref{intro:eq:1}  with special properties focused attention of many authors. The multipeak
solutions of   problem \eqref{intro:eq:1} has been extensively studied among the others by Grossi, Pistoia and Wei
\cite{[GPW]}, Dancer and Yan \cite{[DANCER0]}-\cite{[DANCER1]},  Wang  \cite{[WANG1]} and Yan \cite{[YAN]}.

A multiplicity of solutions of problem \eqref{intro:eq:1} in the presence of symmetries of a compact Lie group has
been studied among the others by Byeon  \cite{[BYEON]}, Vanella \cite{[VANNELLA]}, Wang
\cite{[WANG]}-\cite{[WANG2]}.

The aim of this article is to study connected sets of solutions of problem \eqref{intro:eq:1}.

The first goal of this article is to prove the sufficient conditions for the existence  of solutions of problem
\eqref{intro:eq:1}.

Let $\si = \{0=\lambda_1<\lambda_2<\ldots\}$ denote the set of eigenvalues of the following eigenvalue problem
\begin{equation}\label{intro:eq:2}
\left\{\begin{array}{rcll} -\Delta u  &=& \lambda u& \mathrm{in} \;\o, \\
\ds \frac{\partial u}{\partial \nu} &=& 0 & \mathrm{on} \;
\partial \o,
\end{array}\right.
\end{equation}
and let $\VS(\lambda_i)$ be the eigenspace of  the Laplace operator $-\Delta$ corresponding to the eigenvalue
$\lambda_i \in \si.$

We assume that $f$ is asymptotically linear i.e. $f(x)=f'(\infty)x + o(\p x\p)$, as $\p x\p\ra\infty$ and that $Z=
f^{-1}(0)$ is finite.

In our theorems we put assumptions on $f'(z),$ where $z \in Z \cup \{\infty\}.$   We emphasize that we also treat
problems with resonance at constant solutions and at infinity i.e. it can happen that $f'(z) \in \si$ for some $z
\in Z \cup \{\infty\}.$

The second goal   of this article is to prove the sufficient conditions for continuation  of solutions of the
following problem

\begin{equation}\label{intro:eq:3}
\left\{\begin{array}{rcll} -\Delta u  &=& f(u,\lambda)& \mathrm{in} \;\o, \\
\ds \frac{\partial u}{\partial \nu} &=& 0 & \mathrm{on} \;
\partial \o,
\end{array}\right.
\end{equation}
where $\o\subset \bR^n$ is  an open bounded domain with $C^{1_-}$-boundary    and $f\in C^1(\bR \times \bR,\bR)$.

The third goal of this paper is to study global bifurcations from infinity   of solutions of problem
\eqref{intro:eq:3}.

It is worth in pointing out that application of classical invariants like the Conley index technique and the Morse
theory does not ensure the existence of closed connected sets of critical points of variational problems, see
Ambrosetti \cite{[AMB]}, B\"{o}hme \cite{[BOHME]}, Ize \cite{[IZE0]}, Marino \cite{[MARINO]}, Takens \cite{[TAKENS]}
for examples and discussion.

In other words one can not apply these invariants in order to prove continuation and global bifurcation of solutions
of problem \eqref{intro:eq:3}.

Since the gradient $\nabla \Phi \in C^1(\h1,\h1)$ is of the form compact perturbation of the identity, we apply the
Leray-Schauder degree and the degree for $\sone$-equivariant gradient maps  to the study of critical points
(critical $\sone$-orbits) of the functional $\Phi.$

The choice of the Leray-Schauder degree and the degree for $\sone$-equivariant gradient maps seems to be the best
adapted to our theory.

After this introduction our article is organized as follows.

Since the degree for $\sone$-equivariant gradient maps is not widely known, in Section \ref{sec:prelim} we have
summarized without proofs the relevant material on this invariant, thus making our exposition as self-contained as
possible.

In Section \ref{sec:wlasne} we have studied problem \eqref{intro:eq:2}. In Lemma \ref{lemat:lin:1} we have derived a
formula for the Leray-Schauder degree of the gradient $\nabla_u \Psi \in C^1(\h1 \times \bR,\h1)$ of a functional
$\Psi \in C^2(\h1 \times \bR,\bR)$ associated with problem \eqref{intro:eq:2}. Suppose now that $\bR^n$ is an
orthogonal $\sone$-representation and that $\o \subset \bR^n$ is $\sone$-invariant. Under these assumptions $\h1$ is
an orthogonal $\sone$-representation, the functional $\Psi$ is $\sone$-invariant and its gradient $\nabla_u \Psi$ is
$\sone$-equivariant. In Lemma \ref{lemat:lin:2} we have proved a formula for the degree for $\sone$-equivariant
gradient maps of $\nabla_u \Psi.$

\newpage
In Section \ref{results} our main results are stated and proved.

Subsection \ref{sec:existence} is devoted to the study of the existence of nonconstant solutions of problem
\eqref{intro:eq:1}. In Theorems \ref{tw:ls:1}-\ref{tw:sone:3} we consider non-degenerate case i.e. we assume that
$f'(z) \notin \si$ for every $z \in Z \cup \{\infty\}.$ These theorems ensure the existence of at least one
nonconstant solution of problem \eqref{intro:eq:1}. Notice that in Theorems \ref{tw:sone:1}-\ref{tw:sone:3} we have
assumed that domain $\o$ is $\sone$-invariant.

We emphasize that in the proofs of Theorems \ref{tw:sone:1}-\ref{tw:sone:3} the degree for $\sone$-equi\-va\-riant
gradient maps can not be replaced with the Leray-Schauder degree, see Remark \ref{lsgd}. Additionally, in Theorem
\ref{tw:sone:4} we have proved the existence of at least one nonconstant solution of problem \eqref{intro:eq:1} in a
degenerate case.

In Subsection \ref{sec:continuations} we have studied continuation of nonconstant solutions of problem
\eqref{intro:eq:3}. In Theorems \ref{con:tw:1}, \ref{con:tw:2}, \ref{con:tw:3} we have formulated sufficient
conditions for the existence of closed connected sets of solutions of problem \eqref{intro:eq:3} emanating from a
fixed level $\lambda \in \bR.$

In Subsection \ref{subsec:bif} we have studied global bifurcations from infinity of nonconstant solutions of problem
\eqref{intro:eq:3}. Theorems \ref{bif:tw:1}, \ref{bif:tw:2} are the main theorems of this section.

In Section \ref{examples} we   illustrate the main results of this paper. Namely,    we consider problem
\eqref{intro:eq:1} with $\o=B^2$ and $\o=(0,1)\times B^2.$

\section{Preliminaries}
\label{sec:prelim}

In this section, for the convenience of the reader, we remind the main properties of the degree for $\sone$-equivariant gradient maps defined in
\cite{[RYB1]}. This degree will be denoted briefly by $\nabla_{\sone}\mathrm{-deg}.$

\nt Denote by $\Upsilon(\sone)$ the set of closed subgroups of the group $\sone$ i.e. $\Upsilon(\sone) = \{\sone,
\bZ_1,\bZ_2,\ldots, \bZ_k, \ldots\}.$

\nt Put $\ds U(\sone)= \bZ \oplus \left( \bigoplus_{k=1}^{\infty} \bZ \right)$ and define actions
$$+ , \star :
U(\sone) \times U(\sone) \rightarrow U(\sone),$$
$$\cdot :
\bZ \times U(\sone) \rightarrow U(\sone),$$ as follows

\begin{align}
  \label{doda} \alpha + \beta =&\left(\alpha_0 + \beta_0, \alpha_1 + \beta_1,
 \ldots,\alpha_k+\beta_k,\ldots\right),   \\
 \label{m} \alpha \star \beta =&(\alpha_0 \cdot \beta_0, \alpha_0 \cdot \beta_1 +
\beta_0 \cdot \alpha_1,   \ldots, \alpha_0 \cdot \beta_k + \beta_0 \cdot \alpha_k, \ldots), \\
\gamma \cdot \alpha = &(\gamma \cdot \alpha_0, \gamma \cdot \alpha_1,   \ldots, \gamma \cdot \alpha_k, \ldots),
\end{align}
where $\alpha = (\alpha_0, \alpha_1, \ldots, \alpha_k, \ldots), \beta = (\beta_0, \beta_1, \ldots,\beta_k, \ldots)
\in U(\sone)$ and $\gamma \in \bZ.$ It is easy to check that $(U(\sone),+,\star)$ is a commutative ring with the
unit $\bI=(1,0,\ldots) \in U(\sone)$ and  the trivial element $\Theta=(0,0,\ldots) \in U(\sone).$

The ring $(U(\sone),+,\star)$ is called the Euler ring of the group $\sone$.

\br \label{inv} Notice that $\alpha = (\alpha_0, \alpha_1, \ldots, \alpha_k, \ldots) \in U(\sone)$ is invertible iff
$\alpha_0=\pm 1.$ \er

For a definition of the Euler ring $U(G),$ where $G$ is any compact Lie group, we refer the reader to
\cite{[DIECK]}.

\nt If $\delta_1,\ldots,\delta_q \in U(\sone),$ then we write $\ds
\prod_{j=1}^q \delta_j$ for $\delta_1 \star \ldots \star
\delta_q.$ Moreover, it is understood that $\ds \prod_{j \in
\emptyset} \delta_j= \bI \in U(\sone).$

\nt  Let $\bV$ be a real, finite-dimensional and orthogonal $\sone$-representation. If $v \in \bV,$ then the
subgroup $SO(2)_v=\{g \in \sone : g \cdot v =v\}$ is said to be the isotropy group of $v \in \bV.$ Let $\o \subset
\bV$ be an open, bounded and $\sone$-invariant subset and let $H \in \Upsilon(\sone)$. Then we define
\begin{itemize}
  \item $\o^H=\{v \in \o : H \subset \sone_v\}=\{v \in \o : g v = v \: \forall \: g \in H\},$
  \item $\o_H=\{v \in \o : H = \sone_v\}.$
\end{itemize}

\nt Fix $k \in \bN$ and set
\begin{itemize}
  \item $C^k_{\sone}(\bV,\bR) = \{f \in C^k(\bV,\bR) : f \text{ is } \sone\text{-invariant}\},$
  \item $C^{k-1}_{\sone}(\bV,\bV) = \{f \in C^{k-1}(\bV,\bV) : f \text{ is } \sone\text{-equivariant}\}.$
\end{itemize}
Let $f \in C^1_{\sone}(\bV,\bR).$ Since $\bV$ is an orthogonal $\sone$-representation,  the gradient  $\nabla f \in
C^0_{\sone}(\bV,\bV).$   If $H \in \Upsilon(\sone)$ is a closed subgroup, then $\bV^H$ is a finite-dimensional
$\sone$-re\-pre\-sen\-ta\-tion    and $\big(\nabla f\big)^H =\nabla \big(f_{\mid \bV^H }\big) : \bV^H \rightarrow
\bV^H$ is well-defined $\sone$-equi\-va\-riant gradient map. Choose an open, bounded and $\sone$-invariant subset
$\o \subset \bV$ such that $(\nabla f)^{-1}(0) \cap
\partial \o = \emptyset.$ Under these assumptions we have defined in \cite{[RYB1]} the degree for
$\sone$-equivariant gradient maps $\nabla_{\sone}\mathrm{-deg}(\nabla f,\o) \in U(\sone)$ with coordinates
$$\nabla_{\sone}\mathrm{-deg}(\nabla f,\o)=$$$$=(\nabla_{\sone}\mathrm{-deg}_{\sone}(\nabla f,\o),
\nabla_{\sone}\mathrm{-deg}_{\bZ_1}(\nabla f,\o), \ldots, \nabla_{\sone}\mathrm{-deg}_{\bZ_k}(\nabla f,\o), \ldots
).$$

\br \label{goodrem1} To define the degree for $\sone$-equivariant gradient maps of $\nabla f_0$ we choose (in a
homotopy class of the $\sone$-equivariant gradient map $\nabla f_0$) a  sufficiently good  $\sone$-equivariant
gradient map $\nabla f_1$ and define  this degree for $\nabla f_1.$ The definition does not depend on the choice of
the map $\nabla f_1.$ Roughly speaking the main steps of the definition of the degree for $\sone$-equi\-va\-riant
gradient maps of $\nabla f_0 : (cl(\o),\partial \o) \rightarrow (\bV, \bV \setminus \{0\})$ are the following:
\begin{enumerate}
\item[Step 1.] There is a potential  $f \in C^1_{\sone}(\bV \times [0,1],\bR)$ such that
\begin{enumerate}
\item[(a1)] $(\nabla_v f)^{-1}(0) \cap (\partial \o \times [0,1])=\emptyset,$
\item[(a2)] $\nabla_v f(\cdot,0)=\nabla f_0(\cdot),$
\item[(a3)] $\nabla_v f_1 \in C^1_{\sone}(\bV,\bV),$ where we abbreviate $\nabla_v f(\cdot,1)$ to $\nabla_v f_1,$
\item[(a4)] $(\nabla_v f_1)^{-1}(0) \cap \o^{\sone}=\{v_1,\ldots,v_p\}$ and
\begin{enumerate}
\item[(i)] $\det \nabla^2_{vv} f_1(v_j) \neq 0,$ for all $j =1,\ldots,p,$
\item[(ii)] $\nabla^2_{vv} f_1(v_j)=\left[\begin{array}{cc}
  \nabla^2_{vv}\big(f_1^{\sone}\big)(v_j) & 0 \\
  0 & Id
\end{array} \right] :   \begin{array}{c}
  \bV^{\sone} \\ \oplus \\
  (\bV^{\sone})^{\perp}
\end{array} \longrightarrow  \begin{array}{c}
  \bV^{\sone} \\ \oplus \\
  (\bV^{\sone})^{\perp}
\end{array},$ for all $j =1,\ldots,p,$
\end{enumerate}
\item[(a5)] $(\nabla_v f_1)^{-1}(0) \cap (\o \setminus \o^{\sone})=\{\sone w_1, \ldots, \sone w_q\}$
and
\begin{enumerate}
\item[(i)] $\dim \ker \nabla^2_{vv} f_1(w_j)=1,$ for all $j =1,\ldots,q,$
\item[(ii)] $$\nabla^2_{vv} f_1(w_j)= \left[\begin{array}{ccc}
  0 & 0 &  0\\
  0 & Q_j & 0 \\
  0 & 0 & Id
\end{array}\right]   :$$ $$ \begin{array}{c}
  T_{w_j} (SO(2)w_j) \\ \oplus \\ T_{w_j}(\bV_{\sone_{w_j}}) \ominus  T_{w_j} (SO(2)w_j)\\ \oplus \\
 ( T_{w_j}(\bV_{\sone_{w_j}}))^{\perp}
\end{array}   \longrightarrow  \begin{array}{c}
  T_{w_j} (SO(2)w_j) \\ \oplus \\ T_{w_j}(\bV_{\sone_{w_j}}) \ominus  T_{w_j} (SO(2)w_j)\\ \oplus \\
 ( T_{w_j}(\bV_{\sone_{w_j}}))^{\perp}
\end{array},$$ for all $j =1,\ldots,q.$
\end{enumerate}
\end{enumerate}

\item[Step 2.] The first coordinate of the degree for $\sone$-equivariant gradient maps is defined by
$\ds \dg_{\sone}(\nabla f_0,\o)=\sum_{j=1}^p \sign \det \nabla^2_{vv}(f_1^{\sone})(v_j).$ In other words since
$\nabla \big(f_1^{\sone} \big) = \big(\nabla f_1\big)^{\sone},$ we obtain $$\dg_{\sone}(\nabla f_0,\o)=
\mathrm{deg}_{\mathrm{B}}((\nabla f_1)^{\sone},\o^{\sone},0),$$ where $\mathrm{deg}_{\mathrm{B}}$ denotes  the
Brouwer degree.
\item[Step 3.] Fix $k \in \bN$ and define $$\dg_{\bZ_k}(\nabla f_0,\o)=
\sum_{\{j\in\{1,\ldots,q\} : \sone_{w_j}=\bZ_k\}} \sign \det Q_j,$$
\end{enumerate}
Notice that since $$\mathrm{deg}_{\mathrm{B}}((\nabla f_1)^{\sone},\o^{\sone},0)=\mathrm{deg}_{\mathrm{B}}(\nabla
f_1,\o,0) \textrm{ and } \mathrm{deg}_{\mathrm{B}}(\nabla f_1,\o,0)=\mathrm{deg}_{\mathrm{B}}(\nabla f_0,\o,0)$$
(see \cite{[RAB]}),
 directly by the Step 2.   we obtain $\dg_{\sone}(\nabla f_0,\o)= \mathrm{deg}_{\mathrm{B}}(\nabla
f_0,\o,0).$ Moreover, immediately from the Step 3. we obtain that if $k \in \bN$ and $\sone_v \neq \bZ_k$ for every
$v \in \o,$ then $\dg_{\bZ_k}(\nabla f_0,\o)=0.$ \er

\nt For $\gamma > 0$ and $v_0 \in \bV^{\sone}$ we put $B_{\gamma}(\bV,v_0) = \{v \in \bV :\ \p v - v_0 \p <
\gamma\}$ and $D_{\gamma}(\bV,v_0) = \{v \in \bV :\ \p v - v_0 \p \leq \gamma\}.$ For simplicity of notation we put
$B_\gamma(\bV)=B_\gamma(\bV,0)$ and  $D_\gamma(\bV)=D_\gamma(\bV,0).$

In the following theorem we formulate the main properties of the degree for $\sone$-equivariant gradient maps.

\bt[\cite{[RYB1]}]\label{wlas} Under the above assumptions the degree for $\sone$-e\-qui\-va\-riant gradient maps
has the following properties\et

\begin{enumerate}
  \item  \label{w1}   if
$\nabla_{\sone}\mathrm{-deg}(\nabla f, \o) \neq \Theta,$ then $(\nabla f)^{-1}(0) \cap \o \neq \emptyset,$
\item  if  $\nabla_{\sone}\mathrm{-deg}_H(\nabla f, \o) \neq 0,$ then $(\nabla f)^{-1}(0) \cap \o^H
\neq \emptyset,$
  \item  \label{w3} if
$\o = \o_0 \cup \o_1$ and $\o_0 \cap \o_1 = \emptyset,$ then
$$\nabla_{\sone}\mathrm{-deg}(\nabla f, \o) = \nabla_{\sone}\mathrm{-deg}(\nabla f, \o_0) +
\nabla_{\sone}\mathrm{-deg}(\nabla f, \o_1),$$
  \item  \label{w4}  if
$ \o_0 \subset \o$ is an open $\sone$-invariant subset and $(\nabla f)^{-1}(0) \cap \o \subset \o_0,$ then
$$ \nabla_{\sone}\mathrm{-deg}(\nabla f, \o) = \nabla_{\sone}\mathrm{-deg}(\nabla f, \o_0),$$
\item   \label{w2} if    $ f \in C^1_{\sone}(\bV \times [0,1],\bR)$ is such that
$(\nabla_v f)^{-1}(0) \cap \left(\partial \o \times [0,1] \right) = \emptyset,$ then
$$\nabla_{\sone}\mathrm{-deg}(\nabla f_0, \o)=\nabla_{\sone}\mathrm{-deg}(\nabla f_1,\o),$$

\item \label{wzaw}  if $W$ is an orthogonal $\sone$-representation,  then
  $$\nabla_{\sone}\mathrm{-deg}((\nabla f, Id), \o \times B_{\gamma}(W)) =
\nabla_{\sone}\mathrm{-deg}(\nabla f, \o),$$
\item  if  $ f \in C^2_{\sone}(\bV,\bR)$ is such that $\nabla f(0) = 0$ and $\nabla^2f(0)$
is an   $\sone$-equivariant self-adjoint isomorphism,  then there is $\gamma > 0$ such that
$$\nabla_{\sone}\mathrm{-deg}(\nabla f, B_{\gamma}(\bV)) = \nabla_{\sone}\mathrm{-deg}(\nabla^2 f(0), B_{\gamma}(\bV)).$$
\end{enumerate}

\br \label{goodrem} Directly from the definition of the degree for $\sone$-equivariant gradient maps (see
\cite{[RYB1]}) it follows that
\begin{enumerate}
\item  if $H \in \Upsilon(\sone)$ is a closed subgroup and  $\sone_v \neq H,$ for every $v \in \o,$  then
$\dg_H(\nabla f,\o)=0.$
\item $\dg_{\sone}(\nabla f,\o)= \mathrm{deg}_{\mathrm{B}}(\nabla f, \o,0),$ where $\mathrm{deg}_{\mathrm{B}}$ is  the Brouwer degree.
\end{enumerate}
\er

\nt Below we formulate product formula for the degree for $\sone$-equivariant gradient maps.

\bt[\cite{[RYB2]}] \label{pft} Let $\o_i \subset \bV_i$ be an open, bounded and $\sone$-invariant subset of a
finite-dimensional, orthogonal $\sone$-representation $\bV_i,$  for $i=1,2.$ Let $f_i \in C^1_{\sone}(\bV_i,\bR) $
be such that $\big(\nabla f_i \big)^{-1}(0) \cap
\partial \o_i = \emptyset,$  for $i=1,2.$  Then
$$
\nabla_{\sone}\mathrm{-deg}((\nabla f_1, \nabla f_2), \o_1 \times \o_2) = \nabla_{\sone}\mathrm{-deg}(\nabla f_1,
\o_1) \star \nabla_{\sone}\mathrm{-deg}(\nabla f_2, \o_2).
$$
\et

\nt  For $k \in \bN$ define a map $\rho^k : \sone \rightarrow GL(2,\mathbb{R})$ as follows
\[
\rho^k \left(\left[\begin{array}{lr}
\cos \theta &-\sin \theta\\
\sin \theta &\cos \theta
\end{array}\right] \right)= \left[\begin{array}{lr}
\cos (k\cdot\theta)&-\sin (k\cdot\theta)\\
\sin (k\cdot\theta)&\cos (k\cdot\theta)
\end{array}\right]
\qquad 0\le\theta  <  2\cdot\pi.
\]
For $j,k \in \bN$ we denote by $\mathbb{R}[j,k]$ the direct sum of $j$ copies of $(\mathbb{R}^2 ,\rho^k)$, we also
denote by $\mathbb{R}[j,0]$  the trivial $j$-dimensional $\sone$-representation. We say that two
$\sone$-representations $\bV$ and $\bW$ are equivalent  if there exists an $\sone$-equivariant, linear isomorphism
$T : \bV \rightarrow \bW$. The following classic result gives a complete classification (up to equivalence) of
finite-di\-men\-sio\-nal $\sone$-representations  (see \cite{[ADM]}).

\bt[\cite{[ADM]}] \label{tk} If $\bV$ is a finite-dimensional $\sone$-representation,  then there exist finite
sequences $\{j_i\},\, \{k_i\}$ satisfying:\\ $ (*)\qquad k_i\in \{0\}\cup \bN,\quad  j_i\in \bN,\quad  1\le i\le r,
\: k_1 < k_2 < \dots  < k_r $
\\  such that $\bV$ is equivalent to $\ds\bigoplus^r_{i=1} \mathbb{R}[j_i ,k_i]$. Moreover, the
equivalence class of $\bV$ ($\bV\approx\ds\bigoplus^r_{i=1} \mathbb{R}[j_i,k_i]$) is uniquely determined by
$\{k_i\},\, \{j_i\}$ satisfying $(*)$. \et

\nt Notice that if $\bV\approx\ds\bigoplus^r_{i=1} \mathbb{R}[j_i,k_i]$ and $k_1=0,$ then $\bV^{\sone} \approx
\bR[j_1,0].$ An $\sone$-re\-pre\-sen\-ta\-tion $\bV$ is called nontrivial if $\bV^{\sone} \neq \bV.$ Suppose that
$j' \in \bN, k' \in \bN \cup \{0\}$ and $\bV\approx\ds\bigoplus^r_{i=1} \mathbb{R}[j_i,k_i].$ It is understood that
if $\bR[1,k'] \not \subset \bV,$ then $k' \neq k_i$ for $i=1,\ldots,r.$ Moreover, if $k' \in \bN,$ then
$\bV_{\bZ_{k'}} = \emptyset$ is equivalent to $k' \neq \gcd(k_{i_1},\ldots,k_{i_s})$ for every $\{i_1,\ldots,i_s\}
\subset \{1,\ldots,r\}.$

\nt We will denote by $m^-(L)$ the Morse index of a symmetric matrix $L.$

\nt To apply successfully any degree theory we need computational formulas for this invariant. Below we show how to
compute degree for $\sone$-equivariant gradient maps of a linear, self-adjoint, $\sone$-equivariant isomorphism.

\bl[\cite{[RYB1]}] \label{lindeg} If $\bV \approx \bR[j_0,0] \oplus \bR[j_1,k_1] \oplus \ldots \oplus \bR[j_r,k_r],$
 $L : \bV \rightarrow \bV$ is a self-adjoint, $\sone$-equivariant, linear isomorphism and $\gamma > 0,$ then \el
\begin{enumerate}
  \item $L= \diag (L_0,L_1,\ldots,L_r),$
  \item $$\nabla_{\sone}\mathrm{-deg}_H(L,B_{\gamma}(\bV))=
 \begin{cases}
    (-1)^{m^-(L_0)}, & \text{ for } H = \sone, \\
   \ds (-1)^{m^-(L_0)} \cdot  \frac{m^-(L_i)}{2}, & \text{ for } H = \bZ_{k_i}\\
     0, & \text{ for } H \notin \{\sone, \bZ_{k_1}, \ldots, \bZ_{k_r}\},
  \end{cases}
$$
\item in particular, if $L=-Id,$ then $$\nabla_{\sone}\mathrm{-deg}_H(-Id,B_{\gamma}(\bV))=
 \begin{cases}
    (-1)^{j_0}, & \text{ for } H = \sone, \\
   \ds (-1)^{j_0}  \cdot j_i, & \text{ for } H = \bZ_{k_i},\\
     0, & \text{ for } H \notin \{\sone, \bZ_{k_1}, \ldots, \bZ_{k_r}\}.
  \end{cases}
$$
\end{enumerate}

\nt  Let $(\bH, \langle \cdot,\cdot \rangle _{\bH})$ be an infinite-dimensional, separable Hilbert space which is an
or\-tho\-go\-nal $\sone$-representation    and let $C_{\sone}^1(\bH ,\bR)$ denote the set of $\sone$-invariant
$C^1$-functionals. Fix $\Phi \in C_{\sone}^1(\bH ,\bR)$   such that
\begin{equation}
\label{odwzfun} \nabla \Phi(u)=   u -  \nabla \eta(u),
\end{equation}
where  $\nabla  \eta : \bH   \rightarrow \bH$ is an $\sone$-equivariant   compact operator. Let $\cU \subset \bH$ be
an open, bounded and $\sone$-invariant set such that $\left(\nabla \Phi \right)^{-1}(0) \cap \partial \cU =
\emptyset.$ In this situation $\ds \dg(Id - \nabla \eta, \cU) \in U(\sone)$ is well-defined, see \cite{[RYB1]} for
details  and properties of this degree.

\nt Let $L : \bH \rightarrow \bH$ be a linear, bounded, self-adjoint, $\sone$-equivariant operator with spec\-trum
$\sigma(L)=\{\lambda_i\}.$ By $\bV_L(\lambda_i)$ we will denote eigenspace of $L$ corresponding to the
e\-igen\-va\-lue $\lambda_i$ and we put $\mu_L (\lambda_i)=\dim \bV_L(\lambda_i).$ In other words $\mu_L(\lambda_i)$
is the multiplicity of the eigenvalue $\lambda_i.$ Since operator $L$ is linear, bounded, self-adjoint, and
$\sone$-equivariant, $\bV_L(\lambda_i)$ is a finite-dimensional, orthogonal $\sone$-representation.

\nt For $\gamma > 0$ and $v_0 \in \bH^{\sone}$ put $B_{\gamma}(\bH,v_0) = \{v \in \bH :\  \p v - v_0 \p
 < \gamma\}.$
For simplicity of notation $B_{\gamma}(\bH)$ stands for $B_{\gamma}(\bH,0).$

\nt Combining Theorem 4.5 in  \cite{[RYB1]} with Theorem \ref{pft}
we obtain the following theorem.

\bt  \label{dizom}  Under the above assumptions if $1 \notin \sigma(L),$ then $$\dg(Id - L, B_{\gamma}(\bH))=
\prod_{\lambda_i > 1} \dg(-Id, B_{\gamma}(\bV_L(\lambda_i))) \in U(\sone).$$ It is understood that if $\sigma(L)
\cap [1,+\infty)=\emptyset,$ then $$\dg(Id - L, B_{\gamma}(\bH))=\bI \in U(\sone).$$ \et

\nt Below we formulate the continuation theorem for $\sone$-equivariant gradient o\-pe\-ra\-tors of  the form
compact perturbation of the identity. In other words we study continuation of critical orbits of $\sone$-invariant
$C^1$-functionals. The proof   of this theorem is standard, but in the proof we have to replace the Leray-Schauder
degree with the degree for  $\sone$-equivariant gradient operators.

\bt \label{abscont} Let $\Phi \in C_{\sone}^1(\bH \times \bR,\bR)$ be such that $\ds \nabla_u \Phi(u,\lambda)= u -
\nabla_u \eta(u,\lambda),$ where $\nabla  \eta : \bH  \times \bR \rightarrow \bH$ is an $\sone$-equivariant compact
operator. Fix an open, bounded and $\sone$-invariant subset $\cU \subset \bH$ and $\lambda_0 \in \bR$ such that

\begin{enumerate}
  \item $\left(\nabla_u \Phi(\cdot,\lambda_0 )\right)^{-1}(0) \cap \partial \cU = \emptyset,$
  \item $\ds \nabla_{\sone}-\mathrm{deg}(\nabla_u \Phi(\cdot,\lambda_0),\cU) \neq \Theta \in U(\sone).$
\end{enumerate}
Then there exist continua (closed connected sets) $\cC^{\pm} \subset \bH \times \bR,$ with
$$\cC^- \subset \left(\bH \times (-\infty,\lambda_0]\right) \cap \left(\nabla_u \Phi(\cdot,\lambda_0 )\right)^{-1}(0),$$
$$\cC^+ \subset \left(\bH \times [\lambda_0,+\infty)\right) \cap \left(\nabla_u \Phi(\cdot,\lambda_0 )\right)^{-1}(0),$$
and for both $\cC=\cC^{\pm}$ the following statements are valid
\begin{enumerate}
  \item $\cC \cap (\cU \times \{\lambda_0\}) \neq \emptyset,$
  \item either $\cC$ is unbounded or else $\cC \cap \left((\bH \setminus cl(\cU )) \times \{\lambda_0\}\right) \neq \emptyset.$
\end{enumerate}
\et

>From now on let $\Phi\in C^2_{\sone}(\bH\times\bR,\bR)$  be such that $\nabla_u\Phi(u,\lambda) = u -
\nabla_u\eta(u,\lambda)$, where $\nabla\eta:\bH\times \bR\ra\bH$ is an $\sone$-equivariant compact operator. Fix
$\lambda_+ > \lambda_-$ and assume that there exists $\gamma>0$ such that
\begin{equation}\label{eq:bifcond}
(\nabla \Phi_u(\cdot,\lambda_\pm))^{-1}(0)\subset
B_\gamma(\bH)\times \{\lambda_\pm\}=\emptyset.
\end{equation}
\begin{Definition}\label{pre:def:bif}
An element $\bif \in U(\sone)$ defined as follows
$$
\bif = \dg(\nabla_u\Phi(\cdot,\lambda_+),B_\gamma(\bH)) -
\dg(\nabla_u\Phi(\cdot,\lambda_-),B_\gamma(\bH))
$$
is called the bifurcation index at
$(\infty,[\lambda_-,\lambda_+])$.
\end{Definition}
In the following theorems we have formulated sufficient conditions for the existence of an unbounded closed
connected set of critical orbits bifurcating from infinity. Proofs of this theorems can be found in \cite{[FR]}.
\begin{Theorem}\label{tw:bifinf}Take $\Phi$ as above and let
$\lambda_\pm\in\bR$, $\gamma>0$ be such that condition \eqref{eq:bifcond} holds. If $\bif\not = \Theta\in U(\sone)$,
then there exists an unbounded closed connected component $C$ of
$(\nabla_u\Phi)^{-1}(0)\cap(\bH\times[\lambda_-,\lambda_+])$ such that
$C\cap(B_\gamma(\bH)\times\{\lambda_-,\lambda_+\})\not=\emptyset$.
\end{Theorem}

\newpage
Let $\Phi$ satisfy the following additional assumption:
\begin{itemize}
\item[]$\Phi(u,\lambda) = \frac{1}{2}\lg u,u\rg_\bH -
\frac{1}{2}\lg K_\infty(\lambda)u,u\rg_\bH -
\eta_\infty(u,\lambda)$, where
\begin{enumerate}\item $K_\infty(\lambda):\bH\ra\bH$ is a linear
$\sone$-equivariant self-adjoint operator for every $\lambda\in\bR$,\item the mapping $\bH\times\bR\ni
(u,\lambda)\mapsto K_\infty(\lambda)u\in\bH$ is compact,\item $\nabla_u\eta_\infty:\bH\times \bR\ra \bH$ is an
$\sone$-equivariant compact operator such that $\nabla_u\eta_\infty(u,\lambda)=o(\p u\p)$, as $\p u\p\ra \infty$
uniformly on bounded $\lambda$-intervals.
\end{enumerate}
\end{itemize}
For $\lambda\in\bR$ define $\nabla^2_u\Phi(\infty,\lambda) = Id - K_\infty(\lambda)$. Fix arbitrary
$\lambda_0\in\bR$ and assume that $\ker \nabla_u^2\Phi(\infty,\lambda_0)\not =\{0\}$. Choose $\ep>0$, define
$\lambda_\pm = \lambda_0 \pm \ep$ and assume that the following condition is fulfilled \bc $\{ \lambda\in
[\lambda_-,\lambda_+]:\nabla^2_u\Phi(\infty,\lambda)$  is not  an
 isomorphism$\}=\{\lambda_0\}$. \ec

\begin{Definition}\label{pre:def:meets}
We say that an unbounded closed connected set $\mathcal{C}$ meets $(\infty,\lambda_0)$, if for every $\delta,\eta>0$
$$
\mathcal{C} \cap \{(\h1\setminus B_\gamma(\h1))\times [\lambda_0 - \delta,\lambda_0+\delta]\}\not = \emptyset.
$$
\end{Definition}
The following theorem localize points at which closed connected
sets of solutions of equation $\nabla_u\Phi(u,\lambda)=0$ meet
infinity.
\begin{Theorem}\label{tw:bifinfmeets}
Let $\Phi$ be as above. Choose $\ep,\gamma>0$, $\lambda_0,\lambda_\pm\in \bR$ such that the above conditions are
satisfied. If $\bif\not = \Theta\in U(\sone)$, then the statement of Theorem \ref{tw:bifinf} holds true. Moreover,
$\mathcal{C}$ meets $(\infty,\lambda_0)$.
\end{Theorem}
\section{Linear equation}
\label{sec:wlasne}

Throughout this section we assume that   $\o\subset \bR^n$ is a bounded, open set with $C^{1_-}$-boun\-da\-ry.
Consider the following eigenvalue problem
\begin{equation}\label{eq:eigen:classic}
\left\{\begin{array}{rcll}
       \ds  -\Delta u &=& \lambda u & \;\;in \;\o,
        \\
       \ds \frac{\partial u}{\partial\nu} &=& 0 &\;\; on \;\partial\o.
\end{array}\right.
\end{equation}
Denote by  $\si := \{0=\lambda_1<\lambda_2<\ldots\}$ the set of distinct eigenvalues of problem
\eqref{eq:eigen:classic}. Let $\VS(\lambda_i)$   be the eigenspace of $-\Delta$ corresponding to the eigenvalue
$\lambda_i\in\si$. Additionally define
$$
\nu(\lambda) = \left\{\begin{array}{lcl} \ds
\sum_{\lambda_i<\lambda} \dim \VS(\lambda_i) & if  \;\lambda > 0, \\
0 & if \; \lambda\leq 0.
\end{array}\right.
$$

\nt Solutions of problem \eqref{eq:eigen:classic} are in one to one correspondence with  critical points of
functional $\Psi:\h1\times\bR\ra\bR$ defined by
$$
\Psi(u,\lambda) = \frac{1}{2}\int_\o\p \nabla u\p^2 - \lambda
 u^2 dx.
$$

\nt Computing the gradient $\nabla_u \Psi:\h1\times\bR\ra\h1$ we obtain
$$\lg\nabla_u\Psi(u,\lambda),v\rg_\h1 = \int_\o \nabla u\nabla v -
\lambda uv  dx = $$ $$= \int_\o \nabla u\nabla v + uv - uv - \lambda uv dx = \lg u,v\rg_\h1 - (\lambda + 1) \int_\o
uv dx.$$

\nt According to the Riesz theorem there exists linear bounded operator $\cK:\h1\ra\h1$ given by formula $\ds\lg \cK
u,v\rg_\h1 = \int_\o uv dx$.  By definition $\cK$ is self adjoint and  by the imbedding theorems it is compact.
Hence, $\nabla_u \Psi(u,\lambda)  = u - (\lambda+1)\cK u.$

\nt Fix  $\lambda_i\in\si$ and $u_i\in \VS(\lambda_i)$. Thus $\nabla_u\Psi(u_i,\lambda_i)=0$ and consequently
\begin{equation}\label{eq:postac-liniowego:1}
\nabla_u \Psi(u_i,\lambda) = u_i - (\lambda+1)\cK u_i =  u_i -
\frac{\lambda+1}{\lambda_i+1}u_i=\frac{\lambda_i -
\lambda}{\lambda_i+1}u_i.
\end{equation}

\nt  By the spectral theorem for compact, self-adjoint  operators $\ds \h1 =\overline{\bigoplus_{i=1}^\infty
\VS(\lambda_i)}$. Moreover, for every $u\in\h1$ there exists a unique representation $\ds u = \sum_{i=1}^\infty u_i$
such that $u_i\in \VS(\lambda_i)$ for $i\in\bN \cup \{0\}.$ Hence by \eqref{eq:postac-liniowego:1} we  obtain
\begin{equation}\label{eq:postac-liniowego:2}
\nabla_u\Psi(u,\lambda) = u - (1+\lambda)\cK u = \sum_{i=0}^\infty
\left( \frac{\lambda_i-\lambda}{\lambda_i+1}\right)u_i.
\end{equation}

\nt Since $\nabla_u\Psi(\cdot,\lambda)$ is a family of operators of the form compact perturbation of the identity,
one can apply the Leray-Schauder $\dls$ degree to $\nabla_u\Psi(\cdot,\lambda)$.

\nt The standard proof of the following lemma is omitted.
\begin{Lemma}\label{lemat:lin:1}
Fix $\lambda \not \in\si$ and $\gamma>0$.  Then
$$
\dls(\nabla_u\Psi(\cdot,\lambda),B_{\gamma}(\h1),0) = (-1)^{\nu(\lambda)}.
$$
\end{Lemma}

\begin{Remark}
If  $\lambda \in (0,+\infty) \setminus \si,$ then
$$
\dls(\nabla_u\Psi(\cdot,\lambda),B_{\gamma}(\h1),0) = \prod_{\lambda_i<\lambda}
\dls(-Id,B_{\gamma}(\VS(\lambda_i)),0) \in \{\pm 1\}.
$$
If $\lambda < 0,$ then it is understood that $\dls(\nabla_u\Psi(\cdot,\lambda),B_{\gamma}(\h1),0) =1.$
\end{Remark}

\begin{Remark}
Consider $\bV=\bR^n$ as an orthogonal $\sone$-representation  and let $\o \subset \bV$ be $\sone$-invariant. Then
$\h1$ is an orthogonal $\sone$-re\-pre\-sen\-ta\-tion  with an action given by $(gu)(x)=u(gx).$ For every
$\lambda_i\in\si$, $\VS(\lambda_i)$ is an orthogonal finite-dimensional $\sone$-representation. Moreover, since
$\Psi$ is $\sone$-in\-va\-riant, $\nabla_u\Psi$ is $\sone$-equivariant.
\end{Remark}

\begin{Lemma}\label{lemat:lin:2}
Assume that $\o \subset \bV$ is   $\sone$-invariant. Fix $\lambda\not\in\si$ and $\gamma>0$. Then
\begin{equation*}
\begin{split}
&\dg(\nabla_u\Psi(\cdot,\lambda),B_{\gamma}(\h1))= \prod_{\lambda_i<\lambda} \dg(-Id,B_{\gamma}(\VS(\lambda_i))) \in
U(\sone).
\end{split}
\end{equation*}
It is understood that if $\lambda<0,$ then
$$
\dg(\nabla_u\Psi(\cdot,\lambda),B_{\gamma}(\h1)) = \bI \in U(\sone).
$$
\end{Lemma}
\begin{Proof}
>From \eqref{eq:postac-liniowego:2} it follows that $\sigma((1+\lambda)\cK)=\big\{ \frac{\lambda+1}{\lambda_i+1} :
\lambda_i\in\si\big\}$. By assumption $1\not\in\sigma((1+\lambda)\cK)$. Applying Theorem \ref{dizom} we obtain
$$
\dg(\nabla\Psi_u(\cdot,\lambda),B_{\gamma}(\h1)) = \dg(Id-(\lambda+1)\cK,B_{\gamma}(\h1))= $$$$ =
\prod_{\frac{\lambda+1}{\lambda_i+1} > 1} \dg(-Id,B_{\gamma}(\VS(\lambda_i))),
$$
which completes the proof.
\end{Proof}

\section{Results} \label{results}

In this section we formulate and prove the main results of this article.

In the first subsection we formulate the sufficient conditions for the existence of nonconstant solutions of the
following equation
$$\left\{\begin{array}{rcll} -\Delta u  &=& f(u)& in \;\o, \\
\ds \frac{\partial u}{\partial \nu} &=& 0 & on \;
\partial \o.
\end{array}\right. $$

In the second subsection we study continuation of solutions of the following family of equations

$$
\left\{\begin{array}{rcll} -\Delta u  &=& f(u,\lambda)& in \;\o, \\
\ds \frac{\partial u}{\partial \nu} &=& 0 & on \;
\partial \o,
\end{array}\right.
$$

Finally in the third subsection we study global bifurcations from infinity of solutions of above problem.

In the proofs of theorems of this section as the   topological invariants we use the Leray-Schauder degree and the
degree for \sone-equivariant gradient maps.

\numsubsec
\subsection{Existence of nonconstant solutions} \label{sec:existence}
In this section we study weak solutions of the following equation
\begin{equation}\label{eq:1}
\left\{\begin{array}{rcll} -\Delta u  &=& f(u)& in \;\o, \\
\ds \frac{\partial u}{\partial \nu} &=& 0 & on \;
\partial \o,
\end{array}\right.
\end{equation}
where  $\o\subset \bR^n$ is an open, bounded set with $C^{1_-}$-boundary and $f\in C^1(\bR,\bR)$ satisfy the
following assumption
\begin{itemize}
\item[\textbf{(A.1)}] $\p f'(x)\p \leq a+b\p x\p^p$ for some
$a,b>0$, where $\ds 1 < p < \frac{4}{n-2}$, for $n\geq 3$ and $1 < p < \infty$  for $n=1,2$.
\end{itemize}

\nt Set $F:\bR\ra\bR$  a primitive of $f$ i.e. $\ds F(t)=\int_0^t f(s) ds$. Weak solutions of equation \eqref{eq:1} are in one to one
correspondence with critical points of a  functional $\Phi \in C^2(\h1,\bR)$ defined by $\ds \Phi(u)=\frac{1}{2}\int_\o \p\nabla u\p^2dx - \int_\o
F(u) dx.$

\begin{Remark} \label{rema1}
Constant function $z_0\in \h1$ is a critical point of $\Phi$ iff $z_0\in Z=f^{-1}(0)$. Fix $z_0\in Z.$ Since $
\nabla^2\Phi(z_0) = Id - (1+f'(z_0))\cK$ and \eqref{eq:postac-liniowego:2} it follows that $\nabla^2 \Phi(z_0) : \h1
\rightarrow \h1$ is an isomorphism  iff $ f'(z_0)\not \in \si$.
\end{Remark}

\nt Let us put the following additional assumption
\begin{itemize}
\item[\textbf{(A.2)}] there exists limit $f'(\infty) =
\ds \lim_{\p x\p\ra \infty} \frac{f(x)}{x}$.
\end{itemize}

\nt Notice that $\nabla\Phi(u) = \nabla^2 \Phi(\infty)u + o(\p u\p_\h1)= u - (1+f'(\infty)) \cK u + o(\p u\p_\h1)$
as $\p u\p_\h1\ra\infty$.

\nt We treat $\infty$ as a critical point of $\Phi$. We say that $\infty$ is an isolated critical point of $\Phi$ if
$(\nabla\Phi)^{-1}(0)$ is bounded. Assume that all the elements of $Z \cup \{\infty\}$ are isolated critical points
of $\Phi.$  From now on $\gamma_z$ denotes a positive real number such that:
\begin{itemize}
\item[(i)]  if $z\in Z,$ then $(\nabla\Phi)^{-1}(0) \cap D_{\gamma_z}(\h1,z)=\{z\}$,
\item[(ii)] if $z=\infty,$ then $(\nabla\Phi)^{-1}(0)\subset B_{\gamma_\infty}(\h1).$
\end{itemize}
\begin{Lemma}\label{lemat:degs}
Assume that assumption \textrm{\bf(A.1)} is fulfilled, $z_0\in Z$ and $f'(z_0)\not\in\si$. Then
$\dls(\nabla\Phi,B_{\gamma_{z_0}}(\h1,z_0),0) = (-1)^{\nu(f'(z_0))}.$
\end{Lemma}
\begin{Proof}
It is easy to see that $ \nabla^2\Phi(z_0) = Id - (1+f'(z_0))\cK$. Since $f'(z_0)\not\in\si$, $\nabla^2\Phi(z_0)$ is
an isomorphism. From the properties of the Leray-Schauder degree we get
$\dls(\nabla\Phi,B_{\gamma_{z_0}}(\h1,z_0),0) = \dls(\nabla^2\Phi(z_0),B_{\gamma_{z_0}}(\h1),0)).$ The rest of the
prove is a direct consequence of Lemma \ref{lemat:lin:1}.
\end{Proof}
Since the proof of the next lemma is similar to the proof of Lemma
\ref{lemat:degs},  we will omit it.
\begin{Lemma}\label{lemma:infty:1}
Assume that assumptions \textrm{\bf(A.1)}, \textrm{\bf(A.2)} are satisfied and that
$f'(\infty)\not\in\si$. Then
$\dls(\nabla\Phi,B_{\gamma_\infty}(\h1),0) = (-1)^{\nu(f'(\infty))}.$
\end{Lemma}

\nt Put the following assumptions:
\begin{itemize}
\item[\textbf{(A.3)}] $\#Z<\infty$,
\item[\textbf{(A.4)}] $f'(z)\not\in\si$ for every $z\in Z\cup\{\infty\}$.
\end{itemize}

\nt Define $ Z_+ := \{z\in Z : f'(z)>0\},\;\;Z_-:=\{z\in Z  :  f'(z)<0\}.$

\nt Notice that if assumption (A.4) is fulfilled, then $Z_+ \cup Z_- = Z$.

\nt In the next theorem we ensure the existence of nonconstant solutions of equation \eqref{eq:1}.

\bt \label{tw:ls:1} Suppose that assumptions \textrm{\bf(A.1)}-\textrm{\bf(A.4)} are fulfilled. Moreover, assume that
\begin{enumerate}
\item if $f'(\infty)<0,$ then there exists $z_{0}\in Z_+$ such that $\nu(f'(z_{0})$) is even,
\item if $f'(\infty)>0$ and
$\nu(f'(\infty))$ is odd, then there exists $z_{0}\in Z_+$ such that $\nu(f'(z_{0})$) is even,
\item if $f'(\infty)>0$ and
$\nu(f'(\infty))$ is even, then $\# \{z\in Z_+ : \nu(f'(z))$ is even $\}\not = 1$.
\end{enumerate}
Then there exists at least one nonconstant solution of equation \eqref{eq:1}. \et
\begin{Proof}
By the properties of the Leray-Schauder degree we obtain that
$$\dls(\nabla\Phi,B_{\gamma_\infty}(\h1)\setminus \bigcup_{z\in Z} D_{\gamma_z}(\h1,z),0)=$$
$$=\dls(\nabla\Phi,B_{\gamma_\infty}(\h1),0) - \sum_{z\in Z} \dls(\nabla\Phi,B_{\gamma_z}(\h1,z),0).$$
What is left is to show that
\begin{equation}\nn
\dls(\nabla\Phi,B_{\gamma_\infty}(\h1),0) \not = \sum_{z\in Z} \dls(\nabla\Phi,B_{\gamma_z}(\h1,z),0).
\end{equation}
Suppose, contrary to our claim, that
\begin{equation}\label{eq:proof:4}
\dls(\nabla\Phi,B_{\gamma_\infty}(\h1),0) = \sum_{z\in Z} \dls(\nabla\Phi,B_{\gamma_z}(\h1,z),0).
\end{equation}
By Lemma \ref{lemat:degs} we obtain
\begin{equation}\label{eq:proof:5}
\sum_{z\in Z_-} \dls(\nabla\Phi,B_{\gamma_z}(\h1,z),0) = \#Z_-.
\end{equation}
Now put $Z_+^o:=\{z\in Z_+ : \nu(f'(z))$ is odd$\}$ and $Z_+^e:=\{z\in Z_+ : \nu(f'(z))$ is even$\}$. Then $Z_+ =
Z_+^o\cup Z_+^e$ and $Z_+^o\cap Z_+^e = \emptyset$. Again from Lemma \ref{lemat:degs} it follows that
\begin{equation}\label{eq:proof:6}
\sum_{z\in Z_+} \dls(\nabla\Phi,B_{\gamma_z}(\h1,z),0) = \#Z_+^e - \#Z_+^o.
\end{equation}

By Lemma \ref{lemma:infty:1} we have   if $f'(\infty)<0,$ then $\dls(\nabla\Phi,B_{\gamma_\infty}(\h1),0) = 1$.
Moreover, if  $f'(\infty)>0$ and $\nu(f'(\infty))$ is odd, then $\dls(\nabla\Phi,B_{\gamma_\infty}(\h1),0) = -1$.

\nt Let assumption  (1) or (2) be fulfilled. Then $ \dls(\nabla\Phi,B_{\gamma_\infty}(\h1),0) = -\sign f'(\infty). $ From this and equations
\eqref{eq:proof:4}-\eqref{eq:proof:6} we obtain $\# Z_+^e - \#Z_+^o + \#Z_- = - \sign f'(\infty).$ Moreover, it is easy to see that $\#Z_+ - \#Z_-
= \sign f'(\infty).$ Hence
\begin{equation*}\nn
\left\{\begin{array}{lcr}\# Z_+^e - \#Z_+^o + \#Z_- &=& - \sign f'(\infty), \\
\#Z_+ - \#Z_- &=& \sign f'(\infty). \end{array}\right.
\end{equation*}
We thus get $\#Z_+^e = 0$, a contradiction.

\nt (3)   By Lemma \ref{lemma:infty:1} we obtain $\dls(\nabla\Phi,B_{\gamma_\infty}(\h1),0) = 1.$ Therefore
\begin{equation*}\nn
\left\{\begin{array}{lcr}\# Z_+^e - \#Z_+^o + \#Z_- &=& 1, \\
\#Z_+ - \#Z_- &=& 1, \end{array}\right.
\end{equation*}
which implies $ \#Z_+^e = 1$, a contradiction.
\end{Proof}

\nt From now on we assume that
\begin{itemize}
\item[\textbf{(A.5)}] $\bV=\bR^n$ is a nontrivial orthogonal $\sone$-representation and that $\o \subset \bV$ is
$\sone$-invariant.
\end{itemize}

\br Since $\o \subset \bV$ is $\sone$-invariant, $\h1$ is an orthogonal $\sone$-re\-pre\-sen\-ta\-tion, with
$\sone$-action defined by $(gu)(x)=u(gx),$ and $\Phi \in C^2_{\sone}(\h1,\bR).$   Hence $\nabla \Phi \in
C^1_{\sone}(\h1,\h1).$ \er

\nt The following two  lemmas  are  similar to \ref{lemat:degs},  \ref{lemma:infty:1}, respectively. Since $\nabla \Phi$ is $\sone$-invariant,
instead of the Leray-Schauder degree we will apply the degree for $\sone$-equi\-va\-riant gradient maps.
\begin{Lemma}\label{lemma:deg:sone}
Assume that assumptions \textrm{\bf(A.1)}, \textrm{\bf(A.5)} are fulfilled. Fix $z_0\in Z$ such that $f'(z_0)\not\in\si$. If $z_0\in Z_+,$ then
$$
\dg(\nabla \Phi,B_{\gamma_{z_0}}(\h1,z_0)) = $$
$$=
\prod_{\lambda_i<f'(z_0)}
\dg(-Id,B_{\gamma_{z_0}}(\VS(\lambda_i))) \in U(\sone).
$$
Moreover, if $z_0\in Z_-,$ then $\dg(\nabla\Phi,B_{\gamma_{z_0}}(\h1,z_0)) = \bI \in U(\sone).$
\end{Lemma}
\begin{Proof}
Since $z_0\in\h1$ is a constant function, $B_{\gamma_{z_0}}(\h1,z_0)\subset \h1$ is $\sone$-in\-va\-riant. Moreover,
$\nabla\Phi$ is an $\sone$-equivariant operator of the form compact perturbation of the identity. Hence
$\dg(\nabla\Phi,B_{\gamma}(\h1,z_0))\in U(\sone)$ is well-defined. It is clear that $\nabla^2\Phi(z_0) = Id -
(f'(z_0)+1)\cK$ and that $\nabla^2\Phi(z_0)$ is an isomorphism. From Theorem \ref{wlas} we have
$$
\dg(\nabla\Phi,B_{\gamma_{z_0}}(\h1,z_0)) = \dg(\nabla^2\Phi(z_0),B_{\gamma_{z_0}}(\h1)).
$$
The rest of the proof is a direct consequence of Lemma \ref{lemat:lin:2}.
\end{Proof}

\begin{Lemma}\label{lemma:inf:sone}
Assume that assumptions \textrm{\bf(A.1)}, \textrm{\bf(A.2)}, \textrm{\bf(A.5)}  are satisfied and that
$f'(\infty)\not\in\si$. Then
\begin{enumerate}
  \item if $f'(\infty)>0,$ then
$$
\dg(\nabla\Phi,B_{\gamma_\infty}(\h1)) =$$ $$ \prod_{\lambda_i<f'(\infty)}
\dg(-Id,B_{\gamma_{\infty}}(\VS(\lambda_i))) \in U(\sone),
$$
\item if $f'(\infty)<0,$ then $\dg(\nabla\Phi,B_{\gamma_{\infty}}(\h1)) = \bI \in U(\sone).$
\end{enumerate}
\end{Lemma}
\begin{Proof}
Since  $\nabla\Phi$ is an $\sone$-equivariant operator of the form compact perturbation of the identity and
$\nabla^2\Phi(\infty) = Id - (1+f'(\infty))\cK$ is an isomorphism,
$$\dg(\nabla\Phi,B_{\gamma_{\infty}}(\h1)) = \dg(\nabla^2\Phi(\infty),B_{\gamma_{\infty}}(\h1)).$$ The rest of the proof
is a direct consequence of Lemma \ref{lemat:lin:2}.
\end{Proof}

\nt The following corollary is an immediate consequence of Lemmas  \ref{lindeg}, \ref{lemma:deg:sone},
\ref{lemma:inf:sone}.
\begin{Corollary}\label{nondeg:corollary:1}
If $z \in Z$ and assumptions of Lemma \ref{lemma:deg:sone} are satisfied,  then
\begin{enumerate}
\item if $H \in \Upsilon(\sone)$ and $\dg_H (\nabla\Phi,B_{\gamma_z}(\h1,z))\not = 0,$ then
$$\sign(\dg_H (\nabla\Phi,B_{\gamma_z}(\h1,z))) = (-1)^{ \nu(f'(z))},$$
\item $\dg_{\sone} (\nabla\Phi,B_{\gamma_z}(\h1,z)) = (-1)^{ \nu(f'(z))}.$
\end{enumerate}
If   $z=\infty$ and assum\-ptions of Lemma \ref{lemma:inf:sone} are fulfilled, then
\begin{enumerate}
\item if $H \in \Upsilon(\sone)$ and $\dg_H (\nabla\Phi,B_{\gamma_{\infty}}(\h1))\not = 0,$ then
$$\sign(\dg_H (\nabla\Phi,B_{\gamma_{\infty}}(\h1))) = (-1)^{ \nu(f'(\infty))},$$
\item $\dg_{\sone} (\nabla\Phi,B_{\gamma_{\infty}}(\h1)) = (-1)^{\nu(f'(\infty))}.$
\end{enumerate}
\end{Corollary}

\nt  Define $ \lambda_0 = \min\{\lambda_i\in\si : \VS(\lambda_i) \textrm{ is a nontrivial }
\sone\textrm{-representation}\}.$ Moreover, for $z \in Z \cup \{\infty\}$  define $\ds
\bV(f'(z))=\bigoplus_{\lambda_i < f'(z)} \bV_{-\Delta}(\lambda_i).$

\nt In the next three theorems we prove  the  existence of nonconstant solutions of equation \eqref{eq:1}. Since $\o
\subset \bV$ is $\sone$-invariant, $\nabla \Phi$ is $\sone$-equivariant. Therefore we use in the proofs the degree
for $\sone$-equivariant gradient maps.

It  is worth to point out that we obtain the existence of nonconstant solutions of equation \eqref{eq:1} also if the
assumptions of Theorem \ref{tw:ls:1} are not fulfilled.

\bt \label{tw:sone:1} Suppose that assumptions \textrm{\bf(A.1)}-\textrm{\bf(A.5)} are fulfilled. Moreover, assume that $f'(\infty)<0$ and that
there exists $z_{0}\in Z_+$ such that $\minev < f'(z_{0}).$ Then there exists at least one nonconstant solution of equation \eqref{eq:1}. \et
\begin{Proof}
In view of Theorem \ref{tw:ls:1}, to complete the proof, it is enough to assume that $\nu(f'(z))$ is odd for all
$z\in Z_+$. By the properties of the degree for $\sone$-equivariant gradient maps we obtain
$$\dg(\nabla
\Phi,B_{\gamma_\infty}(\h1)\setminus \bigcup_{z\in Z} D_{\gamma_z}(\h1,z))=$$
$$=\dg(\nabla\Phi,B_{\gamma_\infty}(\h1)) -  \sum_{z\in Z}\dg(\nabla \Phi,B_{\gamma_z}(\h1,z)).$$

\nt Therefore, to complete the proof, it remains to prove that
$$
\dg(\nabla \Phi,B_{\gamma_\infty}(\h1))\not = \sum_{z\in Z}\dg(\nabla \Phi,B_{\gamma_z}(\h1,z)).
$$
Suppose, contrary to our claim, that
\begin{equation} \label{rown}
\dg(\nabla \Phi,B_{\gamma_\infty}(\h1)) = \sum_{z\in Z}\dg(\nabla \Phi,B_{\gamma_z}(\h1,z)).
\end{equation}
Since $\VS(\minev)$ is a nontrivial $\sone$-representation, there is   $k' \in \bN$ such that $\VS(\minev)=$ $=
\bR[1,k']\oplus\bR[1,k']^\bot$. From \eqref{rown} we get
\begin{equation}\label{eq:proof:9}
\dg_{\bZ_{k'}}(\nabla \Phi,B_{\gamma_\infty}(\h1)) = \sum_{z\in Z}\dg_{\bZ_{k'}}(\nabla \Phi,B_{\gamma_z}(\h1,z)).
\end{equation}
Since $f'(\infty)<0$ and  Lemma \ref{lemma:inf:sone}, we obtain
\begin{equation}\label{eq:proof:10}
\dg_{\bZ_{k'}}(\nabla \Phi,B_{\gamma_\infty}(\h1)) = 0.
\end{equation}
If $z\in Z_-,$ then,  by Lemma \ref{lemma:deg:sone}, we have
\begin{equation}\label{eq:proof:11}
\dg_{\bZ_{k'}}(\nabla \Phi,B_{\gamma_z}(\h1,z)) = 0.
\end{equation}
Taking into account  \eqref{eq:proof:9}, \eqref{eq:proof:10} and  \eqref{eq:proof:11}    we obtain
\begin{equation}\label{eq:proof:12}
\sum_{z\in Z_+}\dg_{\bZ_{k'}}(\nabla \Phi,B_{\gamma_z}(\h1,z)) =0.
\end{equation}
Fix $z\in Z_+$. From Lemma \ref{lemma:deg:sone} we have
$$\dg(\nabla \Phi,B_{\gamma_z}(\h1,z))=  \prod_{\lambda_i<f'(z)} \dg(-Id,B_{\gamma_z}(\VS(\lambda_i))) =$$
$$= \dg(-Id,B_{\gamma_z}(\bV(f'(z))).$$

\nt By assumption $\nu(f'(z))$ is odd. Hence from Corollary \ref{nondeg:corollary:1} we obtain

\nt $ \dg_{\sone}(\nabla \Phi,B_{\gamma_z}(\h1,z)) = -1 $ and $ \dg_{\bZ_{k'}}(\nabla \Phi,B_{\gamma_z}(\h1,z)) \leq
0. $ Using the above and \eqref{eq:proof:12} we get $\dg_{\bZ_{k'}}(\nabla\Phi,B_{\gamma}(\h1,z)) = 0 $ for all
$z\in Z_+$. By the assumption there exists $z_0\in Z_+$ such that $f'(z_{0})
> \minev$. Therefore  $\bV(f'(z_0))=\bR[1,k']\oplus \bR[1,k']^\perp$.
Finally, by Lemmas \ref{lindeg}, \ref{lemma:deg:sone},  we obtain
\begin{eqnarray}\nn
\dg_{\bZ_{k'}}(\nabla\Phi,B_{\gamma_{z_0}}(\h1,z_0)) = \dg_{\bZ_{k'}}(-Id,B_{\gamma_{z_0}}(\bV(f'(z_0)))) \not = 0,
\end{eqnarray}
a contradiction.
\end{Proof}

\bt\label{tw:sone:2} Suppose that assumptions \textrm{\bf(A.1)}-\textrm{\bf(A.5)} are fulfilled,  $f'(\infty)>0$ and $\nu(f'(\infty))$ is odd.
Additionally, assume that one of the following conditions is satisfied
\begin{enumerate}
\item there are $z_{0},z_{1}\in Z_+$ such that $f'(z_{0})\geq f'(z_{1})> \minev$ and $f'(z_{0})>f'(\infty)$,
\item there exists exactly one $z_{0}\in Z_+$ such that
\begin{enumerate}
\item $f'(z_{0})>\minev$,
\item there exists $\lambda_{i_0}\in\si$ such that $f'(z_{0})<\lambda_{i_0}<f'(\infty)$ (or
$f'(\infty)<\lambda_{i_0}<f'(z_{0})$) and that $\VS(\lambda_{i_0})$ is a nontrivial $\sone$-representation,
\end{enumerate}
\item there exists $\lambda_{i_0}\in \si$ such that
\begin{enumerate}
\item $f'(z)<\lambda_{i_0}<f'(\infty)$ for all $z\in Z_+$,
\item there exists  $k' \in \bN$ such that
\begin{enumerate}
\item $\VS(\lambda_{i_0}) = \bR[1,k'] \oplus \bR[1,k']^{\perp},$
\item $\bR[1,k'] \not \subset  \VS(\lambda_i)$ for
$\lambda_i\in\si \cap (-\infty,\lambda_{i_0}).$
\end{enumerate}
\end{enumerate}
\end{enumerate}
Then  there exists at least one nonconstant solution of equation \eqref{eq:1}. \et
\begin{Proof}
The proof is similar to that of Theorem \ref{tw:sone:1}. By the properties of the degree for $\sone$-equivariant
gradient maps we obtain
$$\dg(\nabla
\Phi,B_{\gamma_\infty}(\h1)\setminus \bigcup_{z\in Z} D_{\gamma_z}(\h1,z))=$$
$$=\dg(\nabla\Phi,B_{\gamma_\infty}(\h1)) -  \sum_{z\in Z}\dg(\nabla \Phi,B_{\gamma_z}(\h1,z)).$$

\nt It remains to prove that
$$
\dg(\nabla \Phi,B_{\gamma_\infty}(\h1))\not = \sum_{z\in Z}\dg(\nabla \Phi,B_{\gamma_z}(\h1,z)).
$$
Suppose, contrary to our claim, that
\begin{equation}    \nn
\dg(\nabla \Phi,B_{\gamma_\infty}(\h1)) = \sum_{z\in Z}\dg(\nabla \Phi,B_{\gamma_z}(\h1,z)).
\end{equation}
If $z\in Z_-$ and $k \in \bN,$  then, by Lemma \ref{lemma:deg:sone}, we get $ \dg_{\bZ_k}(\nabla
\Phi,B_{\gamma_z}(\h1,z)) = 0 $ and
\begin{equation}\label{eq:proof:13}
\dg_{\bZ_k}(\nabla \Phi,B_{\gamma_\infty}(\h1)) = \sum_{z\in Z_+}\dg_{\bZ_k}(\nabla \Phi,B_{\gamma_z}(\h1,z)).
\end{equation}

\nt From Theorem \ref{tw:ls:1} it follows that, to complete the proof, it suffices to consider the case $\nu(f'(z))$
is odd for all $z\in Z_+\cup \{\infty\}$. Therefore, by Corollary \ref{nondeg:corollary:1}, we obtain that
$\dg_{\bZ_k}(\nabla \Phi,B_{\gamma_z}(\h1,z))  \leq 0$ for all $z\in Z_+\cup\{\infty\}$ and $k \in \bN.$

\nt \textrm{(1)}  Since $\VS(\minev)$ is a nontrivial $\sone$-representation there is  $k' \in \bN$ such that
$\VS(\minev) = \bR[1,k']\oplus \bR[1,k']^\bot.$ Hence, by Lemma \ref{lindeg},  we have
$$
\dg_{\bZ_{k'}}(\nabla \Phi,B_{\gamma_{z_0}}(\h1,z_{0})), \dg_{\bZ_{k'}}(\nabla \Phi,B_{\gamma_{z_1}}(\h1,z_{1}))<0.
$$
Since $f'(z_{0})>f'(\infty)$, it follows that $\bV(f'(\infty)) \subset \bV(f'(z_0))$ and consequently, by Lemmas
\ref{lindeg}, \ref{lemma:deg:sone}, \ref{lemma:inf:sone}, we obtain
$$\dg_{\bZ_{k'}}(\nabla \Phi,B_{\gamma_{z_0}}(\h1,z_{0})) =
\dg_{\bZ_{k'}}(\nabla^2 \Phi(z_0),B_{\gamma_{z_0}}(\h1,z_{0})) =$$
$$= \dg_{\bZ_{k'}}\big(-Id,B_{\gamma_{z_0}} \big(\bV(f'(z_0))\big)\big)  \leq $$
$$\leq
\dg_{\bZ_{k'}}\big(-Id,B_{\gamma_{\infty}} \big(\bV(f'(\infty))\big)\big) =$$ $$ = \dg_{\bZ_{k'}}(\nabla^2
\Phi(\infty),B_{\gamma_{\infty}}(\h1))=  \dg_{\bZ_{k'}}(\nabla \Phi,B_{\gamma_\infty}(\h1)).
$$
Taking together the above inequalities and \eqref{eq:proof:13} we obtain
$$\dg_{\bZ_{k'}}(\nabla \Phi,B_{\gamma_\infty}(\h1)) \geq \dg_{\bZ_{k'}}(\nabla \Phi,B_{\gamma_{z_0}}(\h1,z_{0})) > $$
$$> \dg_{\bZ_{k'}}(\nabla
\Phi,B_{\gamma_{z_1}}(\h1,z_{1}))+\dg_{\bZ_{k'}}(\nabla \Phi,B_{\gamma_{z_0}}(\h1,z_{0}) \geq $$
$$\geq  \sum_{z\in
Z_+}\dg_{\bZ_{k'}}(\nabla \Phi,B_{\gamma_z}(\h1,z)) =\dg_{\bZ_{k'}}(\nabla \Phi,B_{\gamma_\infty}(\h1)),$$ a
contradiction.

\nt \textrm{(2)} Since $\VS(\lambda_{i_0})$ is a nontrivial $\sone$-representation, there is \ $k' \in \bN$ such
that $\VS(\lambda_{i_0})=\bR[1,k']\oplus \bR[1,k']^\bot$. Fix $z \in Z\setminus \{z_0\}.$ Since $f'(z) < \lambda_0,$
$\bV(f'(z))$ is a trivial $\sone$-re\-pre\-sen\-ta\-tion, applying
 Lemmas \ref{lindeg}, \ref{lemma:deg:sone}, we obtain $$\dg_{\bZ_{k'}}(\nabla \Phi,B_{\gamma_z}(\h1,z))=0.$$
Thus
$$\dg_{\bZ_{k'}}(\nabla \Phi,B_{\gamma_\infty}(\h1))  = \sum_{z\in Z} \dg_{\bZ_{k'}}(\nabla
\Phi,B_{\gamma_z}(\h1,z))=$$
$$= \sum_{z\in Z_+} \dg_{\bZ_{k'}}(\nabla
\Phi,B_{\gamma_z}(\h1,z))=\dg_{\bZ_{k'}}(\nabla \Phi,B_{\gamma_{z_0}}(\h1,z_{0})).$$

\nt Let $j_0,j_\infty \in \bN$ be the largest integers such that
$$
\bV(f'(z_0))= \bR[j_0,k']\oplus\bR[j_0,k']^\bot,
 \textrm{ and } \bV(f'(\infty)) = \bR[j_\infty,k']\oplus\bR[j_\infty,k']^\bot.
$$
Since $f'(z_0)<\lambda_{i_0}<f'(\infty)$, we obtain $j_0 < j_\infty.$ Finally, by Lemmas \ref{lindeg},
\ref{lemma:deg:sone}, we obtain
$$
\dg_{\bZ_{k'}}(\nabla \Phi,B_{\gamma_\infty}(\h1)) \not = \dg_{\bZ_{k'}}(\nabla \Phi,B_{\gamma_{z_0}}(\h1,z_{0})),
$$
a contradiction. The same proof remains valid if $f'(\infty)<\lambda_{i_0}<f'(z_{0}).$

\nt \textrm{(3)}  Since  $\ds \bR[1,k'] \not \subset \bV(f'(z))$ for every $z \in Z_+$ and Lemmas \ref{lindeg},
\ref{lemma:deg:sone},
$$
\dg_{\bZ_{k'}}(\nabla \Phi,B_{\gamma_z}(\h1,z)) = \dg_{\bZ_{k'}}(-Id,B_{\gamma_z}(\bV(f'(z))) = 0,
$$
for every $z \in Z.$

\nt Thus, by the above and \eqref{eq:proof:13}, we obtain
$$
\dg_{\bZ_{k'}}(\nabla \Phi,B_{\gamma_\infty}(\h1)) = \sum_{z\in Z_+} \dg_{\bZ_{k'}}(\nabla \Phi,B_{\gamma_z}(\h1,z))
= 0.
$$
Since $\ds \VS(\lambda_{i_0})\subset  \bV(f'(\infty)),$ $\bR[1,k']\subset \bV_\infty$ and consequently
$$
\dg_{\bZ_{k'}}(\nabla \Phi,B_{\gamma_\infty}(\h1))=\dg_{\bZ_{k'}}(-Id,B_{\gamma_\infty}(\bV(f'(\infty)))\not = 0,
$$ a contradiction.
\end{Proof}

\bt\label{tw:sone:3} Suppose that assumptions \textrm{\bf(A.1)}-\textrm{\bf(A.5)} are fulfilled,  $f'(\infty)>0$ and $\nu(f'(\infty))$ is even.
Additionally assume that there exists $z_{0}\in Z_+$ such that $\nu(f'(z_{0}))$ is even and one of the following conditions is fulfilled
\begin{enumerate}
\item there exist $z_{1},z_{2}\in (Z_+ \cup \{\infty\})\setminus
\{z_{0}\}$ such that $f'(z_{1})\geq f'(z_{2})> \lambda_0$ and $f'(z_{1})>f'(z_{0})$, \item there exists exactly one
$z_{1}\in (Z_+ \cup \{\infty\})\setminus \{z_{0}\}$ such that
\begin{enumerate}\item $f'(z_{1})>\minev$, \item there exists $\lambda_{i_0}\in\si$
such that $f'(z_{1})<\lambda_{i_0}<f'(z_{0})$ (or $f'(z_{0})<\lambda_{i_0}<f'(z_{1}))$ and that $\VS(\lambda_{i_0})$
is a nontrivial $\sone$-representation,
\end{enumerate}
\item there exists $\lambda_{i_0}\in \si$ such that
\begin{enumerate} \item $f'(z)<\lambda_{i_0}<f'(z_{0})$ for all $z\in
(Z_+\cup \{\infty\})\setminus\{z_0\}$, \item there exists $k' \in \bN$ such that \begin{enumerate}\item
$\VS(\lambda_{i_0})= \bR[1,k']\oplus \bR[1,k']^\bot$   \item $\bR[1,k'] \not \subset \VS(\lambda_i)$ for
$\lambda_i\in\si \cap (-\infty,\lambda_{i_0}).$
\end{enumerate}\end{enumerate}
\end{enumerate}
Then there   exists at least one nonconstant solution of equation \eqref{eq:1}. \et
\begin{Proof}
The proof is similar to that of Theorem \ref{tw:sone:1}. By the properties of the degree for $\sone$-equivariant
gradient maps we obtain
$$\dg(\nabla
\Phi,B_{\gamma_\infty}(\h1)\setminus \bigcup_{z\in Z} D_{\gamma_z}(\h1,z))=$$
$$=\dg(\nabla\Phi,B_{\gamma_\infty}(\h1)) -  \sum_{z\in Z}\dg(\nabla \Phi,B_{\gamma_z}(\h1,z)).$$

\nt It remains to prove that
$$ \dg(\nabla \Phi,B_{\gamma_\infty}(\h1))\not = \sum_{z\in Z}\dg(\nabla \Phi,B_{\gamma_z}(\h1,z)).$$
Suppose, contrary to our claim, that
\begin{equation}    \nn
\dg(\nabla \Phi,B_{\gamma_\infty}(\h1)) = \sum_{z\in Z}\dg(\nabla \Phi,B_{\gamma_z}(\h1,z)).
\end{equation}
If $z\in Z_-$ and $k \in \bN,$  then, by Lemma \ref{lemma:deg:sone}, we get $ \dg_{\bZ_k}(\nabla
\Phi,B_{\gamma_z}(\h1,z)) = 0 $ and $ \ds \dg_{\bZ_k}(\nabla \Phi,B_{\gamma_\infty}(\h1)) = \sum_{z\in
Z_+}\dg_{\bZ_k}(\nabla \Phi,B_{\gamma_z}(\h1,z)), $ which is equivalent to
\begin{equation}\label{eq:proof:44}
\begin{split}
-\dg_{\bZ_k}(\nabla \Phi,B_{\gamma_{z_0}}(\h1,z_{0})) = \\ =\sum_{z \in Z_+\setminus \{z_{0}\}}\dg_{\bZ_k}(\nabla
\Phi,B_{\gamma_z}(\h1,z))   - \dg_{\bZ_k}(\nabla \Phi,B_{\gamma_\infty}(\h1)).
\end{split}
\end{equation}
 Since $\nu(f'(\infty)),
\nu(f'(z_0))$ are even and Corollary \ref{nondeg:corollary:1}, we
have
$$
- \dg_{\bZ_{k'}}(\nabla \Phi,B_{\gamma_\infty}(\h1))\leq 0, -\dg_{\bZ_{k'}}(\nabla
\Phi,B_{\gamma_{z_0}}(\h1,z_{0}))\leq 0.
$$

\nt Notice that, in view of Theorem \ref{tw:ls:1}, to complete the proof it is enough to consider the case
\begin{equation}\label{opop}
\{z\in Z_+ : \nu(f'(z)) \textrm{ is even} \}  = \{z_0\}.
\end{equation}

\nt \textrm{(1)} Let $z_1,z_2\not=\infty$. Since $\VS(\minev)$ is a nontrivial $\sone$-representation, there is $k'
\in \bN$ such that $\VS(\minev) = \bR[1,k']\oplus \bR[1,k']^\bot.$ Taking into account that $f'(z_1), f'(z_2) >
\lambda_0,\nu(f'(z_1)), \nu(f'(z_1))$ are odd and Corollary \ref{nondeg:corollary:1} we obtain
$$
\dg_{\bZ_{k'}}(\nabla \Phi,B_{\gamma_{z_1}}(\h1,z_{1})), \dg_{\bZ_{k'}}(\nabla \Phi,B_{\gamma_{z_2}}(\h1,z_{2}))<0.
$$
Since $f'(z_{1})>f'(z_0)$, it follows that $\bV(f'(z_0))\subset \bV(f'(z_1))$ and consequently by Lemmas
\ref{lindeg}, \ref{lemma:deg:sone}, \ref{lemma:inf:sone} we obtain
$$\dg_{\bZ_{k'}}(\nabla \Phi,B_{\gamma_{z_1}}(\h1,z_{1})) =
\dg_{\bZ_{k'}}(\nabla^2 \Phi(z_1),B_{\gamma_{z_1}}(\h1,z_{1})) =$$
$$= \dg_{\bZ_{k'}}\big(-Id,B_{\gamma_{z_1}} \big(\bV(f'(z_1))\big)\big)  \leq
- \dg_{\bZ_{k'}}\big(-Id,B_{\gamma_{z_0}} \big(\bV(f'(z_0))\big)\big) =$$
$$ = - \dg_{\bZ_{k'}}(\nabla^2
\Phi(z_0),B_{\gamma_{z_0}}(\h1,z_0))= - \dg_{\bZ_{k'}}(\nabla \Phi,B_{\gamma_{z_0}}(\h1,z_0)).
$$
Taking into account \eqref{eq:proof:44}, \eqref{opop}, Corollary \ref{nondeg:corollary:1} and the above inequalities
we obtain
$$-\dg_{\bZ_{k'}}(\nabla \Phi,B_{\gamma_{z_0}}(\h1,z_0)) \geq \dg_{\bZ_{k'}}(\nabla \Phi,B_{\gamma_{z_1}}(\h1,z_{1})) > $$
$$> \dg_{\bZ_{k'}}(\nabla
\Phi,B_{\gamma_{z_2}}(\h1,z_{2}))+\dg_{\bZ_{k'}}(\nabla \Phi,B_{\gamma_{z_1}}(\h1,z_{1}) \geq $$
$$\geq  \sum_{z\in
Z\setminus \{z_0\}} \dg_{\bZ_{k'}}(\nabla \Phi,B_{\gamma_z}(\h1,z)) - \dg_{\bZ_k}(\nabla
\Phi,B_{\gamma_\infty}(\h1)) = $$$$= - \dg_{\bZ_{k'}}(\nabla \Phi,B_{\gamma_{z_0}}(\h1,z_0)),$$ a contradiction. The
same proof works for $z_1=\infty$ or $z_2=\infty.$ The details are left to the reader.

\nt \textrm{(2)} Assume that $z_1\not=\infty$. Since $\VS(\lambda_{i_0})$ is a nontrivial $\sone$-representation,
there is \ $k' \in \bN$ such that $\VS(\lambda_{i_0})=\bR[1,k']\oplus \bR[1,k']^\bot$. Fix $z \in (Z\cup
\{\infty\})\setminus \{z_0,z_1\}.$ Since $f'(z) < \lambda_0,$ $\ds \bV(f'(z))$ is a trivial
$\sone$-re\-pre\-sen\-ta\-tion and by Lemmas \ref{lindeg}, \ref{lemma:deg:sone} we have $\dg_{\bZ_{k'}}(\nabla
\Phi,B_{\gamma_z}(\h1,z))=0.$ Thus
$$-\dg_{\bZ_{k'}}(\nabla \Phi,B_{\gamma_{z_0}}(\h1,z_0)) =
 \sum_{z\in Z_+\setminus \{z_0\}} \dg_{\bZ_{k'}}(\nabla
\Phi,B_{\gamma_z}(\h1,z)) - $$ $$-\dg_{\bZ_k}(\nabla \Phi,B_{\gamma_\infty}(\h1)) =\dg_{\bZ_{k'}}(\nabla
\Phi,B_{\gamma_{z_1}}(\h1,z_{1})).$$

\nt Let $j_0,j_1 \in \bN$  be the largest integers such that
$$
\bV(f'(z_0))= \bR[j_0,k']\oplus\bR[j_0,k']^\bot,
 \textrm{ and } \bV(f'(z_1))  = \bR[j_1,k']\oplus\bR[j_1,k']^\bot.
$$
Since $f'(z_1)<\lambda_{i_0}<f'(z_0)$, we obtain $j_1 < j_0.$ Finally, by Lemmas \ref{lindeg}, \ref{lemma:deg:sone},
we obtain
$$
-\dg_{\bZ_{k'}}(\nabla \Phi,B_{\gamma_{z_0}}(\h1,z_0)) \not = \dg_{\bZ_{k'}}(\nabla
\Phi,B_{\gamma_{z_1}}(\h1,z_{1})),
$$
a contradiction. The same proof remains valid if $z_1=\infty$ or $f'(z_0)<\lambda_{i_0}<f'(z_{1})$. The details are
left to the reader.

\nt \textrm{(3)}  Since  $\ds \bR[1,k'] \not \subset \bV(f'(z))$ for every $z \in (Z_+\cup\{\infty\})\setminus
\{z_0\}$ and Lemmas \ref{lindeg}, \ref{lemma:deg:sone}, $\dg_{\bZ_{k'}}(\nabla \Phi,B_{\gamma_z}(\h1,z)) =
\dg_{\bZ_{k'}}(-Id,B_{\gamma_z}(\bV(f'(z))) = 0, $ for every $z \in (Z_+\cup\{\infty\})\setminus \{z_0\}$. Thus, by
the above and \eqref{eq:proof:44}, we obtain
$$
-\dg_{\bZ_{k'}}(\nabla \Phi,B_{\gamma_{z_0}}(\h1,z_0)) = $$$$= \sum_{z\in Z_+\setminus \{z_0\}}
\dg_{\bZ_{k'}}(\nabla \Phi,B_{\gamma_z}(\h1,z)) - \dg_{\bZ_{k'}}(\nabla \Phi,B_{\gamma_\infty}(\h1))= 0.
$$
Since $\ds \VS(\lambda_{i_0}) \subset  \ds \bV(f'(z_0)), \bR[1,k']\subset \bV(f'(z_0))$ and consequently
$$
\dg_{\bZ_{k'}}(\nabla \Phi,B_{\gamma_{z_0}}(\h1,z_0))=\dg_{\bZ_{k'}}(-Id,B_{\gamma_{z_0}}( \bV(f'(z_0))))\not = 0,
$$ a contradiction.
\end{Proof}

\begin{Remark} \label{lsgd}
Notice that in Theorems \ref{tw:sone:1}-\ref{tw:sone:3} the degree for $\sone$-equivariant gradient maps can not be
replaced with the  Leray-Schauder degree, since it vanishes. In fact, under assumptions of these theorems it can
happen that
$$
\dls(\nabla \Phi,B_{\gamma_\infty}(\h1)) - \sum_{z\in Z}\dls(\nabla \Phi,B_{\gamma_z}(\h1,z))  = 0 \in \bZ.
$$
and that
$$
\dg(\nabla \Phi,B_{\gamma_\infty}(\h1)) - \sum_{z\in Z}\dg(\nabla \Phi,B_{\gamma_z}(\h1,z)) \not = \Theta \in
U(\sone).
$$
In other words we obtain the existence of nonconstant solution of equation \eqref{eq:1} in the situation when the Leray-Schauder degree is not
applicable i.e. the assumptions of Theorem \ref{tw:ls:1} are not fulfilled.
\end{Remark}

In the rest of this section we consider a degenerate case i.e. we allow $f'(z_0) \in \si$ for some $z_0 \in Z\cup
\{\infty\}.$ To compute a local index of a degenerate isolated critical point of  $\Phi$ we combine the splitting
lemmas (Lemmas 3.2, 3.3 of \cite{[FRR]}) and the product formula for the degree for $\sone$-equivariant gradient
maps, Theorem \ref{pft}.

\nt The following lemma is a consequence of splitting lemmas of \cite{[FRR]}.

\begin{Lemma}\label{lemma:deg:noniso1}
Assume that assumptions \textrm{\bf(A.1), (A.5)} are satisfied. Fix $z_{0} \in Z\cup \{\infty\}$ such that
$f'(z_{0})\in \si$ and $z_{0}\in \h1$ is an isolated critical point of $\Phi$. Then  there exist $\alpha_0>0$ and
$\varphi\in C^2_{\sone}(\VS(f'(z_{0})),\bR)$ such that $0 \in \VS(f'(z_{0}))$ is an isolated critical point of
$\varphi$ and that
$$\dg(\nabla\Phi,B_{\gamma_{z_0}}(\h1,z_{0})) =$$
$$=\dg(\nabla\varphi,B_{\alpha_0}(\VS(f'(z_{0})))) \star
\prod_{\lambda_i<f'(z_{0})}\dg(-Id,B_{\alpha_0}(\VS(\lambda_i))).$$
\end{Lemma}
\begin{Proof}
Fix $z_0 \in Z$ and remind that $\cK :\h1\ra \h1$ is an $\sone$-equivariant, self-adjoint, compact  operator such
that $\lg \cK u,v\rg_\h1 = \int_\o u(x)v(x)dx$. Set $L = (1+f'(z_0))\cK$, $\bV_0 = \ker Id - L$ and $W_0 =
\bV_0^\bot = \im Id - L $. It is easy to see that $\nabla^2\Phi(z_0) = Id - L$ and $\bV_0 = \VS(f'(z_0))$. Define
$\nabla \eta_0:\h1\ra\h1$, $\nabla \eta_0 = \nabla \Phi - (Id-L)$. Then $\nabla\eta_0$ is a compact,
$\sone$-equivariant operator and $\p\nabla\eta_0(u)\p = o(\p u\p)$ as $\p u\p \ra 0$. Now applying Lemma 3.2 of
\cite{[FRR]} we obtain $\alpha_0 > 0$ and $\varphi\in C^2_{\sone}(\VS(f'(z_{0})),\bR)$ with isolated critical point
at the origin and such that
$$ \dg(\nabla\Phi,B_{\alpha_0}(\h1)) = \dg((\nabla\varphi,(Id-L_0)_{|W_0}),B_{\alpha_0}(\bV_0))\times
B_{\alpha_0}(W_0)).
$$
Finally, combining  Theorem \ref{pft} with a slightly modified version of Lemma \ref{lemat:lin:2} (instead of the
operator $\nabla_u \Psi(\cdot,\lambda)$ it is enough to consider the operator $\nabla^2\Phi(z_0)_{\mid W_0}$),   we
obtain
$$
\dg((\nabla\varphi,(Id-L_0)_{|W_0}),B_{\alpha_0}(\bV_0))\times B_{\alpha_0}(W_0)) = $$
$$ = \dg(\nabla\varphi,B_{\gamma_{z_0}}(\VS(f'(z_{0})))) \star \dg((Id - L_0)_{|W_0},B_{\gamma_{z_0}}(W_0))=
$$
$$= \dg(\nabla\varphi,B_{\gamma_{z_0}}(\VS(f'(z_{0})))) \star
\prod_{\lambda_i<f'(z_{0})}\dg(-Id,B_{\gamma_{z_0}}(\VS(\lambda_i))),$$ which completes the proof.

\nt The same proof remains valid for $z_0 = \infty$ but instead of Lemma 3.2 of \cite{[FRR]} we must use Lemma 3.3
of \cite{[FRR]}. The details are left to the reader.
\end{Proof}

\begin{Corollary}\label{corollary:noniso:1}
Fix $z_0 \in Z\cup \{\infty\}$ satisfying assumptions of Lemma \ref{lemma:deg:noniso1} and $k'\in\bN$. Assume that
\begin{enumerate}
\item $\ds \bR[1,k']  \not\subset   \bV(f'(z_0)),$
\item $\VS(f'(z_0))_{\bZ_{k'}} = \emptyset.$
\end{enumerate}
Then  $\dg_{\bZ_{k'}}(\nabla\Phi,B_{\gamma_{z_0}}(\h1,z_0))=0.$
\end{Corollary}
\begin{Proof}
Take $\alpha_0>0$ and $\nabla\varphi$ as in Lemma \ref{lemma:deg:noniso1}. Then
$$\dg(\nabla\Phi,B_{\gamma_{z_0}}(\h1,z_{0}))= $$
$$= \dg(\nabla\varphi,B_{\alpha_0}(\VS(f'(z_{0})))) \star
\dg(-Id,B_{\alpha_0}(\bV(f'(z_0)))).$$ By   $\mathrm{(1)}$ and Lemma \ref{lindeg} we have $\ds
\dg_{\bZ_{k'}}(-Id,B_{\alpha_0}(\bV(f'(z_0))))=0.$ By $\mathrm{(2)}$ and Remark \ref{goodrem} we obtain
$\dg_{\bZ_{k'}}(\nabla\varphi,B_{\alpha_0}(\VS(f'(z_{0}))))=0.$ The rest of the proof is a direct consequence of
product formula \eqref{m}.
\end{Proof}

\bco\label{lemma:deg:noniso2} Fix $z_0 \in Z\cup \{\infty\}$ satisfying assumptions of Lemma
\ref{lemma:deg:noniso1}. If moreover, $\VS(f'(z_0))^{\sone} = \{0\},$ then
$$
\dg_{\bZ_k}(\nabla\Phi,B_{\gamma_{z_0}}(\h1,z_0))= \dg_{\bZ_k}(-Id,B_{\gamma_{z_0}}(\bV(f'(z_0)))),
$$
for all $k\in\bN$ such that $\VS(f'(z_0))_{\bZ_k}=\emptyset$. \eco
\begin{Proof} Take $\alpha_0>0$ and $\nabla\varphi$ as in Lemma \ref{lemma:deg:noniso1}. Then
$$\dg(\nabla\Phi,B_{\gamma_{z_0}}(\h1,z_{0}))= $$
$$= \dg(\nabla\varphi,B_{\alpha_0}(\VS(f'(z_{0})))) \star
\dg(-Id,B_{\alpha_0}(\bV(f'(z_0)))).$$ Since $\VS(f'(z_0))^{\sone} = \{0\},$
$\dg_{\sone}(\nabla\varphi,B_{\gamma_{z_0}}(\VS(f'(z_0)),0)) = 1.$  Moreover, since $\VS(f'(z_0))_{\bZ_k}=\emptyset$
and   Remark \ref{goodrem}, $\dg_{\bZ_{k}}(\nabla\varphi,B_{\gamma_{z_0}}(\VS(f'(z_0)),0)) = 0.$ The rest of the
proof is a direct consequence of formula \eqref{m}.
\end{Proof}

\nt We can now proof the analog of Theorems \ref{tw:sone:1}, \ref{tw:sone:2}, \ref{tw:sone:3}. It is worth to point
out that  in this theorem we allow $\ds \bigcup_{z \in Z \cup \{\infty\}} \{f'(z)\} \cap \si \neq \emptyset.$

\bt\label{tw:sone:4} Let assumptions \textrm{\bf(A.1)}-\textrm{\bf (A.3)}, \textrm{\bf(A.5)} be fulfilled. Moreover, assume that there are
$z_{0}\in Z\cup \{\infty\}$, $\lambda_{i_0}\in \si$ and $k'\in\bN$ such that
\begin{enumerate}
\item either $f'(z_{0})\not \in \si $ or $\VS(f'(z_0))^{\sone} =
\{0\}$ and $\VS(f'(z_0))_{\bZ_{k'}} = \emptyset$,
\item$f'(z_{0}) > \lambda_{i_0}>f'(z)$ for all $z \in (Z\cup\{\infty\}) \setminus \{z_{0}\},$
\item $\bR[1,k'] \subset \VS(\lambda_{i_0}),$
\item $\VS(\lambda_i)_{\bZ_{k'}}=\emptyset$ for all $\lambda_i\in\si$,
$\lambda_i<\lambda_{i_0}$.
\end{enumerate}
Then there exists at least one nonconstant solution of equation \eqref{eq:1}. \et
\begin{Proof}
Without loss of generality we can assume that elements of  $Z\cup \{\infty\}$ are isolated critical points of the
potential $\Phi.$ To complete the proof it is enough to show that
\begin{equation}\nn
\dg(\nabla \Phi,B_{\gamma_\infty}(\h1))\not = \sum_{z\in Z}\dg(\nabla \Phi,B_{\gamma_z}(\h1,z)).
\end{equation}
Suppose, contrary to our claim, that
\begin{equation}\label{eq:proof:15}
\dg(\nabla \Phi,B_{\gamma_\infty}(\h1))= \sum_{z\in Z}\dg(\nabla \Phi,B_{\gamma_z}(\h1,z)).
\end{equation}
Fix $z\in (Z\cup \{\infty\})\setminus \{z_{0}\}$. By assumptions $\mathrm{(2),(4)}$ we obtain that
$\VS(\lambda_i)_{Z_{k'}}=\emptyset$ for all $\lambda_i\in\si \cap (-\infty,f'(z)].$ Therefore $\ds \bR[1,k'] \not
\subset \bV(f'(z))$ and if $f'(z)\in\si,$ then $\VS(f'(z))_{\bZ_{k'}} = \emptyset$. Hence, from Corollary
\ref{corollary:noniso:1}, we obtain
\begin{equation}\label{eq:proof:14}
\dg_{Z_{k'}}(\nabla \Phi,B_{\gamma_z}(\h1,z)) = 0,
\end{equation}
for all $z\in  (Z\cup \{\infty\})\setminus \{z_{0}\}$.

\nt If $f'(z_0) \not \in \si,$ then, by Lemmas \ref{lemma:deg:sone} ($z_0 \in Z$) or Lemma \ref{lemma:inf:sone}
($z_0=\infty$), we get
$$
\dg_{Z_{k'}}(\nabla \Phi,B_{\gamma_{z_0}}(\h1,z_{0})) = \dg_{Z_{k'}}(-Id,B_{\gamma_{z_0}}(\bV(f'(z_0)),z_0)).
$$

\nt If $f'(z_0) \in \si,$ then, by  assumption $\mathrm{(1)}$ and Corollary \ref{lemma:deg:noniso2}, we get
$$
\dg_{Z_{k'}}(\nabla \Phi,B_{\gamma_{z_0}}(\h1,z_{0})) = \dg_{Z_{k'}}(-Id,B_{\gamma_{z_0}}(\bV(f'(z_0)),z_0)).
$$
Finally, since $\lambda_{i_0}<f'(z_0)$ and $\bR[1,k'] \subset  \VS(\lambda_{i_0})$, we obtain $\ds \bR[1,k'] \subset
\bV(f'(z_0)).$

\nt Thus, by Lemma \ref{lindeg}, we obtain $\ds \dg_{Z_{k'}}(-Id,B_{\gamma_{z_0}}(\bV(f'(z_0)),z_0)) \not = 0$ and
consequently
\begin{equation}\label{ggt}
\dg_{Z_{k'}}(\nabla \Phi,B_{\gamma_{z_0}}(\h1,z_{0}))   \not= 0.
\end{equation}
Combining   \eqref{eq:proof:14} with \eqref{ggt} we get
$$\dg(\nabla \Phi,B_{\gamma_\infty}(\h1,0)) - \sum_{z\in Z} \dg(\nabla
\Phi,B_{\gamma_z}(\h1,z))= $$
$$=- \dg(\nabla \Phi,B_{\gamma_{z_0}}(\h1,z_{0})) \not =0,$$
contrary to \eqref{eq:proof:15}.
\end{Proof}
\begin{Remark}\label{noniso:remark:1}
Let us notice that in the above theorem the same proof works for assumption $\mathrm{(4)}$ replaced by assumptions
\begin{enumerate}
  \item $\bR[1,k'] \not \subset  \VS(\lambda_i)$ for all $\lambda_i \in \si \cap (-\infty, \lambda_{i_0}),$
  \item $\VS(f'(z))_{\bZ_{k'}}=\emptyset$ for all $z\in (Z\cup \{\infty\}) \cap \si.$
\end{enumerate}
The details are left to the reader.
\end{Remark}

\numsubsec
\subsection{Continuation of   solutions} \label{sec:continuations} Consider a family of
equations of the form
\begin{equation}\label{eq:zparam}
\left\{\begin{array}{rcll} -\Delta u  &=& f(u,\lambda)& in \;\o, \\
\ds \frac{\partial u}{\partial \nu} &=& 0 & on \;
\partial \o,
\end{array}\right.
\end{equation}
where $f \in C^1(\bR \times \bR,\bR), f(\cdot,\lambda)$ satisfies condition \textbf{(A.1)} for every $\lambda\in\bR$
and $\o \subset \bR^n$ is an open, bounded set with $C^{1_-}$-boundary.

In this section we study continuation of nonconstant solutions of family \eqref{eq:zparam}.

\begin{Remark}\label{con:rem:1}
Consider a  functional $\Phi \in C^2(\h1\times\bR,\bR)$ defined as follows
$$
\Phi(u,\lambda) = \int_\o \p\nabla u(x)\p^2 - F(u(x),\lambda) dx,
$$
where $F'_u=f.$

Define $Z_0 = (f(\cdot,0))^{-1}(0)$ and assume that
\begin{enumerate}
\item $\# Z_0 < \infty,$
\item all the elements of $Z_0\cup \{\infty\}$ are isolated critical points of $\Phi(\cdot,0).$
\end{enumerate}
Define an open bounded set $\cU$ in the following way  $$ \cU = B_{\gamma_\infty}(\h1)\setminus
 \bigcup_{z\in Z_0} D_{\gamma_z} (\h1,z).$$ Since  $(\nabla_u \Phi(\cdot,0))^{-1}(0)\cap
\partial \cU = \emptyset, (\nabla_u \Phi(\cdot,0))^{-1}(0)\cap \partial \cU = \emptyset.$
Therefore, by the properties of the Leray-Schauder degree, we obtain
$$
\dls(\nabla_u\Phi(\cdot,0),\cU,0) = $$ $$=\dls(\nabla_u\Phi(\cdot,0),B_{\gamma_\infty}(\h1),0) - \sum_{z\in
Z_0}\dls(\nabla_u\Phi(\cdot,0),B_{\gamma_z}(\h1,z),0).
$$
If moreover, assumption $\mathbf{(A.5)}$ is fulfilled, then $\Phi \in C^2_{\sone}(\h1\times\bR,\bR)$ and $\cU$ is
$\sone$-invariant. Therefore, by the properties of the degree for $\sone$-equivariant gradient maps, we obtain
$$
\dg(\nabla\Phi_u(\cdot,0),\cU) = $$ $$=\dg(\nabla_u\Phi(\cdot,0),B_{\gamma_\infty}(\h1)) -   \sum_{z\in
Z_0}\dg(\nabla_u\Phi(\cdot,0),B_{\gamma_z}(\h1,z)).
$$
\end{Remark}

\bt \label{con:tw:1} Fix $f \in C^1(\bR \times \bR,\bR)$ and
assume that $f(\cdot,0)$   satisfies assumptions  of Theorem
\ref{tw:ls:1}. Then there exist closed connected sets $\cC^{\pm}$
such that
$$\cC^- \subset (\h1 \times (-\infty,0])\cap(\nabla_u\Phi)^{-1}(0) \:  and  \: \cC^+ \subset
(\h1\times [0,+\infty))\cap(\nabla_u\Phi)^{-1}(0).$$ Moreover, for $\cC=\cC^{\pm}$
\begin{itemize}
\item[(i)] $\ds \cC\cap \big(\big(B_{\gamma_\infty}(\h1)\setminus
\bigcup_{z\in Z_0} D_{\gamma_z}(\h1,z)\big)\times
\{0\}\big)\not=\emptyset$, \item[(ii)] either $\cC$ is unbounded
or $\cC\cap (Z_0\times \{0\})\not=\emptyset$.
\end{itemize}
\et
\begin{Proof} Repeating the reasoning from the proof of Theorem \ref{tw:ls:1} we obtain

$$
\dls(\nabla_u\Phi(\cdot,0),B_{\gamma_\infty}(\h1),0) \not = \sum_{z\in Z_0} \dls
(\nabla_u\Phi(\cdot,0),B_{\gamma_z}(\h1,z),0).
$$
Define $\ds \cU = B_{\gamma_\infty}(\h1)\setminus   \bigcup_{z\in Z_0} D_{\gamma_z}(\h1,z)$ and notice that   $
\dls(\nabla_u\Phi(\cdot,0),\cU,0)\not = 0. $ Applying Theorem \ref{abscont} we obtain the existence of closed
connected sets $\cC^\pm$  such that
\begin{eqnarray}\nn
&\cC^-&\subset(\h1\times
(-\infty,0])\cap(\nabla_u\Phi)^{-1}(0),\\\nn &\cC^+&\subset
(\h1\times [0,+\infty))\cap(\nabla_u\Phi)^{-1}(0),
\end{eqnarray}
$\cC=\cC^\pm$ satisfies (i) and    either $\cC$ is unbounded or else $\displaystyle C\cap \left((\h1\setminus
cl(\cU))\times \{0\}\right)\not=\emptyset$. \vspace{3mm}
\\
By definition $\ds \left(\h1\setminus cl(\cU)\right) \cap (\nabla_u\Phi(\cdot,0))^{-1}(0)\subset \bigcup_{z\in Z_0}
B_{\gamma_z}(\h1,z).$ On the other hand $\displaystyle\bigcup_{z\in Z_0} B_{\gamma_z}(\h1,z)\cap
(\nabla_u\Phi(\cdot,0))^{-1}(0)= Z_0,$ which completes the proof.
\end{Proof}

\bt \label{con:tw:2} Fix $f \in C^1(\bR \times \bR,\bR)$ and assume that $f(\cdot,0)$ satisfies assumptions of one
of Theorems \ref{tw:sone:1}, \ref{tw:sone:1}, \ref{tw:sone:3}.   Then there exist closed connected sets $\cC^{\pm}$
such that
$$\cC^- \subset (\h1\times (-\infty,0])\cap(\nabla_u\Phi)^{-1}(0) \: and \: \cC^+ \subset (\h1\times
[0,+\infty))\cap(\nabla_u\Phi)^{-1}(0).$$ Moreover, for $\cC=\cC^{\pm}$
\begin{itemize}
\item[(i)] $\ds \mathcal{C} \cap \big(\big(B_{\gamma_\infty}(\h1)\setminus
\bigcup_{z\in Z_0} D_{\gamma_z}(\h1,z)\big)\times
\{0\}\big)\not=\emptyset$,
\item[(ii)] either $\mathcal{C} $ is unbounded or
$\mathcal{C} \cap (Z_0\times \{0\})\not=\emptyset$.
\end{itemize}
\et
\begin{proof}
Repeating the reasoning from the proofs of Theorems \ref{tw:sone:1}-\ref{tw:sone:3} we obtain
$$
\dg(\nabla_u\Phi(\cdot,0),B_{\gamma_\infty}(\h1)) \not = \sum_{z\in
Z_0}\dg(\nabla_u\Phi(\cdot,0),B_{\gamma_z}(\h1,z)).
$$
Set $\ds \cU = B_{\gamma_\infty}(\h1)\setminus   \bigcup_{z\in Z_0} D_{\gamma_z}(\h1,z)$. Notice that by  Remark
\ref{con:rem:1} we obtain
\\ $\dg(\nabla\Phi(\cdot,0),\cU) \not = \Theta$. The rest of the proof is a direct consequence of Theorem \ref{abscont}.
\end{proof}

\bt \label{con:tw:3} Assume   that $f(\cdot,0)$   satisfies assumptions  of Theorem \ref{tw:sone:4}. Then there
exists infinite sequence of nonconstant solutions of equation \eqref{eq:zparam} with $\lambda=0$ converging to some
$z\in Z_0$  or there exist closed connected sets $\cC^{\pm}$ such that
$$\cC^-\subset(\h1\times (-\infty,0])\cap(\nabla_u\Phi)^{-1}(0) \: and \: \cC^+\subset (\h1\times
[0,+\infty))\cap(\nabla_u\Phi)^{-1}(0).$$ Moreover, for $\cC=\cC^{\pm}$
\begin{itemize}
\item[(i)] $\ds \mathcal{C} \cap \big(\big(B_{\gamma_\infty}(\h1)\setminus
\bigcup_{z\in Z_0} D_{\gamma_z}(\h1,z)\big)\times
\{0\}\big)\not=\emptyset$,
\item[(ii)] either $\mathcal{C}$ is unbounded or
$\mathcal{C} \cap (Z_0\times \{0\})\not=\emptyset$.
\end{itemize}
\et
\begin{Proof}
Suppose that   doesn't exist a sequence of nonconstant solutions of equation \eqref{eq:zparam} with $\lambda=0$
converging to some point in $Z_0$. Then all the points $z\in Z_0$ are isolated critical points of $\Phi(\cdot,0)$.
Repeating the reasoning from the proof of Theorem \ref{tw:sone:4} we obtain
$$
\dg(\nabla_u\Phi(\cdot,0),B_{\gamma_\infty}(\h1)) \not =
\sum_{z\in Z_0}\dg(\nabla_u\Phi(\cdot,0),B_{\gamma_z}(\h1,z)).
$$
We set $\ds \cU =B_{\gamma_\infty}(\h1)\setminus   \bigcup_{z\in Z_0} D_{\gamma_z}(\h1,z).$  Applying  Remark
\ref{con:rem:1} we obtain the following  $\dg(\nabla\Phi,\cU) \not = \Theta$.  The rest of the proof is a direct
consequence of Theorem \ref{abscont}.
\end{Proof}

\numsubsec
\subsection{Bifurcations from infinity}\label{subsec:bif}
In this section we study bifurcations from infinity of solutions of a family of equations of the form

\begin{equation}\label{eq:zparam1}
\left\{\begin{array}{rcll} -\Delta u  &=& f(u,\lambda)& in \;\o, \\
\ds \frac{\partial u}{\partial \nu} &=& 0 & on \;
\partial \o,
\end{array}\right.
\end{equation}
where $f \in C^1(\bR \times \bR,\bR), f(\cdot,\lambda)$ satisfies condition \textbf{(A.1)} for every
$\lambda\in\bR.$ $\o \subset \bR^n$ is an open, bounded set with $C^{1_-}$-boundary. Moreover, we assume that
assumption  \textbf{(A.5)} is fulfilled.

\noindent \textbf{(B.1)} Fix $\lambda_+>\lambda_-$ and assume that $f(\cdot,\lambda_\pm)$ satisfy assumption
\textbf{(A.2)} and $f'(\infty,\lambda_\pm)\not\in\si$.

Notice that under such an assumption $\nabla^2_u\Phi(\infty,\lambda_\pm) = Id - f'(\infty,\lambda_\pm)\cK$, where
operator $\cK:\h1\ra\h1$ is given by the formula $\displaystyle \lg \cK u ,v\rg_\h1 = \int_\o u(x) v(x) dx$.
Moreover, $\nabla^2_u\Phi(\infty,\lambda_\pm)$ is a linear isomorphism iff $f'(\infty,\lambda_\pm)\not\in\si$.
Therefore operator $\nabla_u\Phi(\cdot,\lambda_\pm)$ is asymptotically linear at infinity and its derivative at
infinity is a linear isomorphism. Thus there exists $\gamma>0$ such that
\begin{equation}\label{bif:eq:1}
(\nabla_u\Phi(\cdot,\lambda_\pm))^{-1}(0)\subset
B_\gamma(\h1)\times \{\lambda_\pm\}
\end{equation}
and we can define $\bif \in U(\sone).$ The following theorem is a direct consequence of Theorem \ref{tw:bifinf}.

\newpage
\begin{Theorem}\label{bif:tw:1}
Let $\lambda_+ > \lambda_-$ satisfy assumption \textrm{\bf(B.1)} and fix $\gamma>0$ such that condition
\eqref{bif:eq:1} holds. If $\bif \not = \Theta \in U(\sone)$, then there exists an unbounded closed connected
component $\mathcal{C}$ of $(\nabla_u\Phi)^{-1}(0)\cap (\h1\times [\lambda_-,\lambda_+])$ such that $\mathcal{C}
\cap (B_\gamma(\h1)\times \{\lambda_-,\lambda_+\})\not = \emptyset$.
\end{Theorem}
In the next lemma show how to verify that the  bifurcation index $\bif \in U(\sone)$ is nontrivial.
  In this lemma, for simplicity, we assume that
$f'(\infty,\lambda_+)>f'(\infty,\lambda_-)$. It is clear that similar lemma can be formulated if
$f'(\infty,\lambda_+) < f'(\infty,\lambda_-)$.
\begin{Lemma}\label{bif:lemma:1}
Let $\lambda_+ > \lambda_-$ satisfy assumption \textrm{\bf(B.1)} and $f'(\infty,\lambda_+)>f'(\infty,\lambda_-)$.
Then $\bif\not=\Theta \in U(\sone)$ iff at least one of the following conditions is satisfied:
\begin{enumerate}\item there exists $\lambda_{i_0}\in\si$, such that
$f'(\infty,\lambda_-) < \lambda_{i_0} < f'(\infty,\lambda_+)$ and
\\ $\VS(\lambda_{i_0})\not = \VS(\lambda_{i_0})^{\sone}$, \item
$\displaystyle\sum_{f'(\infty,\lambda_-)< \lambda_i <
f'(\infty,\lambda_+)} \dim\VS(\lambda_i)$  is odd.
\end{enumerate}
\end{Lemma}
\begin{proof}
Fix $\gamma>0$ such that condition \eqref{bif:eq:1} holds.
Directly from the definition we have
$$
\bif = $$$$ = \dg(\nabla_u\Phi(\cdot,\lambda_+),B_\gamma(\bH)) -
\dg(\nabla_u\Phi(\cdot,\lambda_-),B_\gamma(\bH)).
$$
By Lemma \ref{lemma:inf:sone} we obtain
$$
\bif = \left(\prod_{f'(\infty,\lambda_-)< \lambda_i <
f'(\infty,\lambda_+)} \dg(-Id,B_{\gamma_{\infty}}(\VS(\lambda_i)))
- \bI\right)\star $$$$
\star\dg(\nabla_u\Phi(\cdot,\lambda_-),B_\gamma(\bH)).
$$

Since $f'(\infty,\lambda_-)\not \in\si$ and Lemma \ref{lemma:inf:sone},
$\dg(\nabla_u\Phi(\cdot,\lambda_-),B_\gamma(\bH))$ is  the degree of isomorphism.

Thus, by Lemma \ref{dizom}, $\dg_{\sone}(\nabla_u\Phi(\cdot,\lambda_-),B_\gamma(\bH))  = \pm 1.$ Consequently by
Remark \ref{inv} we  obtain that  $\dg(\nabla_u\Phi(\cdot,\lambda_-),B_\gamma(\bH))$ is invertible in $U(\sone)$ and
therefore $\bif = \Theta \in U(\sone)$ iff
$$
\prod_{f'(\infty,\lambda_-)< \lambda_i < f'(\infty,\lambda_+)}
\dg(-Id,B_{\gamma_{\infty}}(\VS(\lambda_i))) =  \bI.
$$
This equality can be rewritten as follows
\begin{equation}\label{ewq}
\dg(-Id,B_{\gamma_{\infty}}(\bigoplus_{f'(\infty,\lambda_-)< \lambda_i < f'(\infty,\lambda_+)}\VS(\lambda_i))) =
\bI.
\end{equation}

By Lemma \ref{dizom} condition \eqref{ewq} is satisfied  iff $\displaystyle \bigoplus_{f'(\infty,\lambda_-)<
\lambda_i < f'(\infty,\lambda_+)}\VS(\lambda_i)$ is a trivial and even dimensional $\sone$-representation.

Therefore, $\bif \not = \Theta \in U(\sone)$ iff $\displaystyle \bigoplus_{f'(\infty,\lambda_-)< \lambda_i <
f'(\infty,\lambda_+)}\VS(\lambda_i)$ is a nontrivial $\sone$-representation or is odd dimensional, which completes
the proof.
\end{proof}
\begin{Corollary}
Fix  $\lambda_\pm\in\bR$ and $\eta>0$ as in Theorem \ref{bif:tw:1}. Additionally, assume that
$f'(\infty,\lambda_+)>f'(\infty,\lambda_-)$ and that one of the following conditions is satisfied:
\begin{enumerate}\item there exists $\lambda_{i_0}\in\si$, such that
$f'(\infty,\lambda_-) < \lambda_{i_0} < f'(\infty,\lambda_+)$ and
\\ $\VS(\lambda_{i_0})\not = \VS(\lambda_{i_0})^{\sone}$, \item
$\displaystyle\sum_{f'(\infty,\lambda_-)< \lambda_i <
f'(\infty,\lambda_+)} \dim\VS(\lambda_i)$  is odd.
\end{enumerate}
Then the statement of Theorem \ref{bif:tw:1} holds true.
\end{Corollary}
\begin{Remark}
It is easy to see that the analogous corollary as above holds true
also if $f'(\infty,\lambda_+) < f'(\infty,\lambda_-)$.
\end{Remark}

We will need the following assumption.

\noindent \textbf{(B.2)} Let $f(\cdot,\lambda)$ satisfy assumption \textbf{(A.2)} for all $\lambda\in\bR$.

Notice that under assumption \textbf{(B.2)} $\nabla_u\Phi(\cdot,\lambda)$ is asymptotically linear for all
$\lambda\in\bR$. Hence $\nabla_u^2\Phi(\infty,\lambda)$ is defined for all $\lambda\in\bR$ and
$\nabla^2_u\Phi(\infty,\lambda) = Id - f'(\infty,\lambda)\cK$. Let $\lambda_0\in\bR$ be such that
$f'(\infty,\lambda_0)\in\si$. Choose $\ep>0$, define $\lambda_\pm=\lambda_{0}\pm\ep$ and assume that
\begin{equation}\label{bif:eq:2}
\{\lambda \in [\lambda_-,\lambda_+] : f'(\infty,\lambda)\in \si\}
= \{\lambda_{0}\}.
\end{equation}

It is clear that under this condition we have \bc $\{\lambda \in [\lambda_-,\lambda_+] :
\nabla^2_u\Phi(\infty,\lambda)$ is not an isomorphism$\} = \{\lambda_{0}\}$. \ec We can consider $\bif \in
U(\sone)$. The following lemma is a direct consequence of Lemma \ref{bif:lemma:1}.
\begin{Lemma}
Let assumption \textrm{\bf (B.2)} holds and let $\lambda_0,\lambda_\pm\in\bR$ be such that condition
\eqref{bif:eq:2} is satisfied. Then $\bif \not = \Theta \in U(\sone)$ iff one of the following conditions  holds:
\begin{enumerate}
\item $\VS(f'(\infty,\lambda_0))\not =
\VS(f'(\infty,\lambda_0))^{\sone}$,\item
$\VS(f'(\infty,\lambda_0))$ is odd.
\end{enumerate}
\end{Lemma}

The next theorem is the consequence of the above lemma  and
Theorem \ref{tw:bifinfmeets}.
\begin{Theorem}\label{bif:tw:2}
Let assumption \textrm{\bf (B.2)} holds and let $\lambda_0,\lambda_\pm\in\bR$ be such that condition \ref{bif:eq:2}
is satisfied. Moreover, assume that $\VS(f'(\infty,\lambda_0))\not = \VS(f'(\infty,\lambda_0))^{\sone}$ or $\dim
\VS(f'(\infty,\lambda_0))$ is odd. Then the statement of Theorem \ref{bif:tw:1} holds true and moreover
$\mathcal{C}$ meets $(\infty,\lambda_{0})$.
\end{Theorem}
\section{Examples} \label{examples}
\numsec

In this section we illustrate the abstract results proved in the previous section.

Define  $\bV_1=\bR[1,1]$, $\bV_2=\bV_1\oplus\bR[1,0]$ and denote by  $\o_1 \subset \bV_1$ an open disc of radius one
in $\bV_1$ and $\o_2=\o_1\times(0,1)\subset \bV_2$. Since $\sone$-representations $\bV_1,\bV_2$ are orthogonal, sets
$\o_1, \o_2$ are $\sone$-invariant. First we remind some standard facts about  $ \sigma(-\Delta,\o_i)$, $i=1,2$.   \\
Throughout this section we assume that $k,n \in \bN \cup \{0\}.$ Moreover, if $k \in \bN,$ then $n \in \bN$ and by
$x_{kn}$ we denote the $n$-th solution of $J_k'(x)=0$ in $(0,+\infty)$, where $J_k$ is an $k$-th Bessel function. If
$k=0,$ then $n \in \bN \cup \{0\}$ and by $x_{0n}$ we denote the $n$-th solution of $J_0'(x)=0$ in $[0,+\infty).$
Notice that $x_{00}=0.$

\begin{Lemma}(\cite{[MATFIZ]}) \label{matfiz}
Under the above assumptions
\begin{enumerate}
\item $\ds \sigma(-\Delta,\o_1) = \{\lambda_{kn}=x_{kn}^2\}_{k=1,n=1}^\infty \cup \{\lambda_{0n}=x_{0n}^2\}_{n=0}^{\infty};$
with corresponding eigenvectors   in spherical coordinates given by
\begin{enumerate}
\item if $k >0,$ then $n> 0$ and  $\ds
\lambda_{kn} \longrightarrow v_{kn}(r,\phi)=J_k(x_{kn} r)\left\{\begin{array}{l} \cos k \phi,\\ \sin k
\phi,\end{array}\right. $
\item if $k=0,$   then $\lambda_{0n} \longrightarrow v_{0n}(r,\phi)=J_0(x_{0n} r).$
\end{enumerate}
\item
$\ds \sigma(-\Delta,\o_2)=\{\lambda_{knj}=(\pi n)^2 + x_{kj}^2 \}_{k=1,n=0,j=1}^\infty \cup \{\lambda_{0nj}=(\pi
n)^2 + x_{0j}^2 \}_{n=0,j=0}^\infty;$ with corresponding eigenvectors  in cylindrical coordinates given by
\begin{enumerate}
\item if $k >0,$ then $j >0$ and $\ds \lambda_{knj}   \longrightarrow v_{knj} (r,\phi,z) = \cos(n\pi z)J_k(x_{kj}r)
\left\{\begin{array}{l} \cos k\phi,\\ \sin k\phi,
\end{array}\right.$
\item if $k=0,$ then $\ds \lambda_{0nj} \longrightarrow v_{0nj} (r,\phi,z) = \cos(n\pi z)J_0(x_{0j}r).$
\end{enumerate}
\end{enumerate}
\end{Lemma}

In the next lemma we show some properties of zeros of derivatives of Bessel functions.

\bl \label{lemma:ill:1} Under the above assumptions
\begin{enumerate}
\item $0=x_{00} < x_{01} < x_{02}< \ldots,$ \item $0 < x_{k1} <
x_{k2}< x_{k3} <\ldots$, for $k \in \bN,$ \item $x_{11} < x_{21} <
x_{31}< \ldots$,
\end{enumerate}
\el

Applying    Lemmas \ref{matfiz}, \ref{lemma:ill:1}   we obtain the following corollary.

\begin{Corollary}\label{besprop}Under the above assumptions
\begin{enumerate}\item $\lambda_{00} < \lambda_{01} < \ldots$,
\item $\lambda_{11}<\lambda_{21}<\ldots$, \item
$\lambda_{101}<\lambda_{201}<\ldots$.
\end{enumerate}
\end{Corollary}

In   Lemma \ref{o1} we describe eigenspaces of $-\Delta$  corresponding to eigenvalues $\sigma(-\Delta;\o_1)$ as
$\sone$-representations.
\begin{Lemma}\label{o1}
If $k\in \bN \cup \{0\}, n \in \bN$ and $\lambda_{kn} \in \sigma(-\Delta;\o_1),$ then $\bR[1,k] \subset
\VS(\lambda_{kn})$. Additionally, $\VS(\lambda_{00})=\bR[1,0].$
\end{Lemma}
\begin{Proof}  First of all notice that from Lemma \ref{matfiz} we obtain
\begin{enumerate}
\item if $k>0,$ then
$\mathrm{span}_{\bR}\{J_k(x_{kn} r)\cos k \phi,J_k(x_{kn} r)\sin k \phi\} \subset \bV_{-\Delta}(\lambda_{kn}),$
\item if $k=0$ and $n > 0,$ then $\mathrm{span}_{\bR}\{J_0(x_{0n} r)\} \subset \bV_{-\Delta}(\lambda_{0n}),$
\item if $k=0$ and $n = 0,$ then $\mathrm{span}_{\bR}\{v_{00}\} = \bV_{-\Delta}(\lambda_{00}).$
\end{enumerate}
Since the $\sone$-action $\sone \times \bH^1(\o_1) \rightarrow \bH^1(\o_1)$ is given by
$$\left(\left[\begin{array}{lr}
\cos \theta &-\sin \theta\\
\sin \theta &\cos \theta
\end{array}\right],u \right)(r,\phi)=u(r,\phi+\theta),$$ it is easy to check that
\begin{enumerate}
  \item $\mathrm{span}_{\bR}\{J_k(x_{kn} r)\cos k \phi,J_k(x_{kn} r)\sin k \phi\} \approx \bR[1,k],$
  \item $\mathrm{span}_{\bR}\{J_0(x_{0n} r)\} \approx \bR[1,0],$
  \item $\mathrm{span}_{\bR}\{v_{00}\} \approx \bR[1,0],$
\end{enumerate}
which completes the proof.
\end{Proof}

The corollary below is a consequence of Lemma \ref{o1}.

\begin{Corollary}\label{corollary:ill:1}
If $\lambda\in\sigma(-\Delta,\o_1)$ and
$$
 \{(k,n)\in (\bN\cup\{0\})^2 : \lambda_{kn}\in\sigma(-\Delta,\o_1)\textrm{ and }
\lambda_{kn}=\lambda\}= \{(k_1,n_1),  \ldots,(k_s,n_s)\},
$$
then $\VS(\lambda) \simeq \bR[1,k_1]\oplus \bR[1,k_2] \oplus \cdots \oplus \bR[1,k_s]$. \\ Moreover, if  $a>0$ and
$\ds \nu(a) = \sum_{\lambda_{kn}<a} \dim \VS(\lambda_{kn}),$ then $$\nu(a) \textrm{ is even  iff  } \#\{\lambda_{0n}
: \lambda_{0n}<a\} \textrm{ is even }.$$
\end{Corollary}

In   Lemma \ref{eigens} we describe eigenspaces of $-\Delta$  corresponding to eigenvalues $\sigma(-\Delta;\o_2)$ as
$\sone$-representations.
\begin{Lemma} \label{eigens}
If $k,n,j \in \bN \cup \{0\}$ and $\lambda_{knj} \in \sigma(-\Delta;\o_2),$ then $\bR[1,k] \subset
\VS(\lambda_{knj})$. Additionally, $\VS(\lambda_{000})=\bR[1,0].$
\end{Lemma}
\begin{proof} In fact the proof is the same as the proof of Lemma \ref{o1}. The details are left to the reader.\end{proof}

The following corollary is a direct consequence of Lemmas  \ref{eigens}.

\begin{Corollary}\label{corollary:ill:2}
If $\lambda\in\sigma(-\Delta,\o_2)$ and
$$
\{(k,n,j)\in (\bN\cup\{0\})^3\} : \lambda_{knj}\in\sigma(-\Delta,\o_2) \textrm{ and } \lambda_{knj}=\lambda\} =
\{(k_1,n_1,j_1), \ldots,(k_s,n_s,j_s)\},
$$
then $\VS(\lambda) \simeq \bR[1,k_1]\oplus \bR[1,k_2] \oplus \cdots \oplus \bR[1,k_s].$ \\ Moreover, if  $a>0$ and
$\ds \nu(a) = \sum_{\lambda_{knj}<a} \dim \VS(\lambda_{knj}),$ then $$\nu(a) \textrm{ is even iff }
\#\{\lambda_{0nj} : \lambda_{0nj}<a\} \textrm{ is even }.$$
\end{Corollary}

\begin{Remark}
For $i=1,2$ let $\minev(\o_i)$  be the smallest eigenvalue in $\sigma(-\Delta,\o_i)$, such that $\VS(\minev(\o_i))$
is a nontrivial $\sone$-representation. It clear that $$\minev(\o_1) = \lambda_{11}=x_{11}^2 =
\lambda_{101}=\minev(\o_2).$$
\end{Remark}

The proof of the   lemma below is a direct consequence of estimations from \cite[Section 15.3, p.486]{[Watson]}.
\begin{Lemma}\label{lemma:ill:2}
For every $k\in\bN$, $\lambda_{k1}\in \sigma(-\Delta,\o_1)$ and
$\lambda_{k01}\in\sigma(-\Delta,\o_2)$ we have
$$
k(k+2)<\lambda_{k1}<2k(k+1),\;\;k(k+2)<\lambda_{k01}<2k(k+1).
$$
Consequently, $3<\minev(\o_1)=\minev(\o_2) < 4$. \end{Lemma}

\begin{Example}
Consider equation
\begin{equation}\label{ill:eq:2}
\left\{\begin{array}{rcll} -\Delta u  &=& f(u)& \textrm{ in } \o_1, \\
\ds \frac{\partial u}{\partial \nu} &=& 0 & \textrm{ on }
\partial \o_1,
\end{array}\right.
\end{equation}
where $f$ satisfies the following assumptions:
\begin{enumerate}
\item $f\in C^1(\bR,\bR)$,
\item $\mid f'(x) \mid \leq a+b\p x\p^q$ for some $a,b>0, q\in \bN$,
\item $f(t) = f'(\infty)t + o(\mid t\mid)$, where $\mid t\mid\ra\infty$,
\item  $50 < f'(\infty) < 99,$
\item $\#Z<\infty,$ where $Z=f^{-1}(0),$ \item $f'(z) \not \in
\sigma(-\Delta,\o_1)$, for all $z\in Z,$
\item \label{seven} there are $z_0,z_1\in Z$ such that $4 <  f'(z_1) < 99 < f'(z_0).$
\end{enumerate}
It is clear that $f$ satisfies assumptions \textbf{(A.1)}-\textbf{(A.5)} of the previous section. Moreover, it is
known that   $\lambda_{02}=x_{02}^2 \approx 49 < f'(\infty) < 100 \approx  x_{03}^2 =\lambda_{03}.$ Therefore by
Corollary \ref{corollary:ill:1} we obtain  that $\nu(f'(\infty))$ is odd.

Taking into account assumption  \eqref{seven} and Lemma \ref{lemma:ill:2} we obtain that $f'(z_0) > f'(z_1)
> \lambda_0(\o_1)$ and $f'(z_0) > f'(\infty).$ Now it is easy to verify that under the above assumptions $f$ satisfies assumption (1) of Theorem
\ref{tw:sone:2}. Thus  there exists at least one nonconstant weak solution of equation \eqref{ill:eq:2}.

If there exists exactly one $z_0\in Z$ such that $f'(z_0) > \lambda_0(\o_1),$ then in order to use Theorem
\ref{tw:sone:2} we have to replace assumption \eqref{seven} with the following assumption:
\begin{itemize}
\item[(7')] there exists $k',n'\in\bN$ such that $f'(z_0) <
\lambda_{k'n'} < 50$ ($99 < \lambda_{k'n'} < f'(z_0)$).
\end{itemize}
Indeed, with assumption \eqref{seven} replaced by assumption (7') the assumption (2) of Theorem \ref{tw:sone:2} is
fulfilled.
\end{Example}

\begin{Example}
Consider equation
\begin{equation}\label{ill:eq:1}
\left\{\begin{array}{rcll} -\Delta u  &=& f(u)& \textrm{ in } \o_2, \\
\ds \frac{\partial u}{\partial \nu} &=& 0 & \textrm{ on }
\partial \o_2,
\end{array}\right.
\end{equation}
where $f$ satisfies the following assumptions:
\begin{enumerate}
\item $f\in C^1(\bR,\bR)$,
\item $\mid f'(x)\mid \leq a+b\p x\p^{3}$ for some $a,b>0$,
\item $f(t) = f'(\infty)t + o(\mid t
\mid)$, where $\mid t\mid\ra\infty$, \item  $f'(\infty)<0$, \item
$\#Z<\infty,$ where $Z=f^{-1}(0),$ \item $f'(z) \not \in
\sigma(-\Delta,\o_2)$, for all $z\in Z$, \item there exists
$z_0\in Z$ such that $f'(z_0) > 4$.
\end{enumerate}
It is clear that $f$ satisfies assumptions \textbf{(A.1)}-\textbf{(A.5)} of the previous section.  Moreover, by
assumption (7) and Lemma \ref{lemma:ill:2}, $f'(z_0)>\minev(\o_2)$. Applying Theorem \ref{tw:sone:1} we obtain
nonconstant weak solutions of equation \eqref{ill:eq:1}.

Define $Z_+ = \{z\in Z \mid \: f'(z) > 0\}$ and assume that $0< f'(z) < 9,$ for every $z\in Z_+.$ Since $x_{01}
\simeq 3.83, f'(z) < \lambda_{001} = x_{01}^2$ and $f'(z) < \lambda_{010} = \pi^2$. By Lemmas \ref{matfiz},
\ref{lemma:ill:1} we obtain $\{\lambda_{0nj} \in \sigma(-\Delta,\o_2) : \lambda_{0nj}<9 \} = \{\lambda_{000}\}$.
Therefore, by Corollary \ref{corollary:ill:1}, we obtain that  $\nu(f'(z))$ is odd for every $z\in Z_+$. Notice that
assumptions of Theorem \ref{tw:ls:1} are not fulfilled. In other words we can not apply the Leray-Schauder degree to
obtain the existence of nonconstant weak solutions of equation \eqref{ill:eq:1}.
\end{Example}

\begin{Example}
Consider equation \eqref{ill:eq:1} and assume that
\begin{enumerate}
\item $f\in C^1(\bR,\bR)$,
\item $\mid f'(x) \mid \leq a+b\p x\p^{3}$ for some $a,b>0$,
\item $f(t) = f'(\infty)t + o(\mid t \mid)$, where $\mid t \mid\ra\infty$,
\item $\# Z<\infty$,
\item there exists $z_0\in Z$ and $k'\in\bN$ such that:
\begin{enumerate}\item $f'(z_0)\not \in
\sigma(-\Delta,\o_2)$, \item $f'(z_0)> 2k'(k'+1)$, \item
$f'(z)<k'(k'+2)$ for $z \in (Z\cup\{\infty\})\setminus \{z_0\}$.
\end{enumerate}
\end{enumerate}
It is clear that $f$ satisfies assumptions \textbf{(A.1)}, \textbf{(A.2)} and \textbf{(A.3)} of the previous
section.  Combining assumptions (5.b), (5.c) with Lemma \ref{lemma:ill:2} we obtain that $f'(z) <  \lambda_{k'01} <
f'(z_0)$ for all $z\in (Z\cup\{\infty\})\setminus\{z_0\}$. Fix $\lambda\in\sigma(-\Delta,\o_1) \cap
(0,\lambda_{k'01}).$ By Corollary \ref{corollary:ill:2}, $\VS(\lambda) \simeq \bR[1,k_1]\oplus \ldots \oplus
\bR[1,k_s]$ for some $k_1,\ldots ,k_s\in \bN \cup \{0\}$.

\nt We claim that $k_i < k'$ for every $1 \leq i \leq s.$ Suppose, contrary to our claim,  that $k_{i_0} \geq k'$
for some $1 \leq i_0 \leq s$. Then, by Corollary \ref{corollary:ill:2}, there exist $n_{i_0}\in\bN\cup \{0\},
j_{i_0} \in \bN$ such that $\lambda_{k_{i_0}n_{i_0}j_{i_0}} = \lambda$. From Lemmas \ref{matfiz}, \ref{lemma:ill:1}
we obtain
$$\lambda_{k_{i_0}n_{i_0}j_{i_0}}=(\pi n_{i_0})^2+x_{k_{i_0}j_{i_0}}^2 \geq x_{k_{i_0}j_{i_0}}^2 \geq x_{k_{i_0}1}^2 =
\lambda_{k_{i_0}01}.$$ By Corollary \ref{besprop} we obtain $\lambda=\lambda_{k_{i_0}n_{i_0}j_{i_0}}\geq
\lambda_{k_{i_0}01} \geq \lambda_{k'01}$, a contradiction. Thus $k_{i} < k'$ for $i=1,\ldots,s$ and consequently
$\bV(\lambda)_{\bZ_{k'}} = \emptyset.$ Taking into account assumption  (5.c) and Lemma \ref{lemma:ill:2} we obtain
$\VS(f'(z))_{\bZ_{k'}}=\emptyset$ for all $z\in (Z\cup\{\infty\})\setminus \{z_0\}$ such that $f'(z)\in
\sigma(-\Delta,\o_2)$.

Notice that all the assumptions of Theorem \ref{tw:sone:4} are satisfied. Applying this theorem we obtain the
existence of at least one nonconstant weak solutions of equation \eqref{ill:eq:1}.

Suppose now that (5.a) does not hold, i.e. $f'(z_0)\in \sigma(-\Delta,\o_2)$. In order to obtain the existence of
weak nonconstant solutions of equation \eqref{ill:eq:1} we have to assume:
\begin{itemize}
\item[(5.a')] $\VS(f'(z_0))^{\sone} = \{0\}$ and $\VS(f'(z_0))_{\bZ_{k'}} = \emptyset$.
\end{itemize}
It is clear that under the above assumption and assumptions (1)-(3), (5.b) and (5.c) Theorem \ref{tw:sone:4} holds.
This assumption is equivalent to the following one
\begin{itemize}
\item[(5.a'')] $f'(z_0)\not = \lambda_{0nj}$ for $n,j\in\bN\cup
\{0\}$ and $f'(z_0) \not = \lambda_{k''nj}$ where $k'' = k' m$ for
$m\in \bN$, $n \in \bN\cup \{0\}$, $j\in\bN$.
\end{itemize}
\end{Example}

\begin{Example}
In this example we illustrate bifurcations from infinity. Consider
the family of equations
\begin{equation}\label{ill:eq:3}
\left\{\begin{array}{rcll} -\Delta u  &=&  f(u,\lambda)& \textrm{ in } \o_1, \\
\ds \frac{\partial u}{\partial \nu} &=& 0 & \textrm{ on }
\partial \o_1,
\end{array}\right.
\end{equation}
where $f$ satisfies the following assumptions:
\begin{enumerate}
\item $f\in C^1(\bR\times\bR,\bR)$, \item $\mid f'(x,\lambda) \mid
\leq a+b\p x\p^{q}$ for some $a,b>0, q\in \bN$ and all
$\lambda\in\bR$, \item there exist limits $f'(\infty,\lambda_\pm)$
for some $\lambda_+,\lambda_- > 0$,\item
$f'(\infty,\lambda_\pm)\not\in\sigma(-\Delta,\o_1)$, \item
$f'(\infty,\lambda_-)<f'(\infty,\lambda_+)$.
\end{enumerate}
Under the above assumptions $\bif \not = \Theta \in U(\sone)$ iff one of the following conditions is satisfied:
\begin{itemize}
\item[(a)] there exists $\lambda_{kn}\in \sigma(-\Delta,\o_1)$,
$n,k\in\bN$ such that $f'(\infty,\lambda_+) > \lambda_{nk} >
f'(\infty,\lambda_-)$, \item[(b)] $\#\{\lambda_{0n}\in
\sigma(-\Delta,\o_1):f'(\infty,\lambda_-)<\lambda_{0n} <
f'(\infty,\lambda_+)\}$ is odd.
\end{itemize}
Therefore, all the assumptions of Theorem \ref{bif:tw:1} are fulfilled if one of the conditions (a), (b) is
satisfied.

Let us now replace condition (3) with the following:
\begin{itemize}
\item[(3')] there exists limit $f'(\infty,\lambda)$,  for all
$\lambda\in\bR$.
\end{itemize}
Let $\lambda_0$ be such that $f'(\infty,\lambda_0)\in\si$. If there exists $k,n\in\bN$ such that
$f'(\infty,\lambda_0) = \lambda_{kn}$, then by Corollary \ref{corollary:ill:1} we get
$\VS(f'(\infty,\lambda_0))\not=\VS(f'(\infty,\lambda_0))^{\sone}$. Otherwise, $f'(\infty,\lambda_0) = \lambda_{0n}$
for some $n\in\bN$ and from Corollaries \ref{besprop} and \ref{corollary:ill:1} we conclude that $\dim
\VS(f'(\infty,\lambda_0)) = 1$. Hence, if there exists $\lambda_-<\lambda_+$ such that
$$
\{\lambda\in[\lambda_-,\lambda_+] : f'(\infty,\lambda)\in\si\} =
\{\lambda_0\},
$$
then $\bif\not = \Theta \in U(\sone)$. Applying  Theorem \ref{bif:tw:2}  we obtain the existence of an unbounded
connected set $\mathcal{C}$ of weak solutions of equation \eqref{ill:eq:3} which  meets $(\infty,\lambda_0)$.

Now let $f(x,\lambda) = \lambda f(u)$ and assume that $f'(\infty)
= f'(\infty,1)\not = 0$. Fix
$\lambda_{i_0}\in\sigma(-\Delta,\o_1)$ and put
$\displaystyle\lambda_0 = \frac{\lambda_{i_0}}{f'(\infty)}$. Then
there exists $\ep
>0$ and $\lambda_\pm = \lambda_0 \pm\ep$ such
$$
\{\lambda\in[\lambda_-,\lambda_+] : \lambda f'(\infty)\in\si\} =
\{\lambda_0\}.
$$
Consequently, $\bif\not = \Theta \in U(\sone)$. By  Theorem \ref{bif:tw:2} there exists an unbounded connected set
$\mathcal{C}$ of  weak solutions of equation \eqref{ill:eq:3} which meets $(\infty,\lambda_0)$.
\end{Example}

\end{document}